\documentclass[A4j,11pt]{article}
\usepackage{amsfonts,amsmath,amscd}
\usepackage{graphics}
\usepackage{amssymb}
\usepackage[all]{xy}
\begin{document}

\title{{\bf Examples of certain kind of\\
minimal orbits of Hemann actions}}
\author{{\bf Naoyuki Koike}}
\date{}
\maketitle

\begin{abstract}
We give examples of certain kind of minimal orbits of Hermann actions and 
discuss whether each of the examples is austere.  
\end{abstract}

\vspace{0.5truecm}

%\vspace{0.2truecm}

%\noindent
%$\overline{\qquad\qquad\qquad\qquad\qquad\qquad}$

%\vspace{0.1truecm}

%{\rm \footnotesize{2000 Mathematics Subject Classification. Primary 53C40, 
%53C35.}}

%{\rm \footnotesize{Key words and phrases. equifocal submanifold, polar 
%action.}}

\section{Introduction}
Let $N=G/K$ be a symmetrc space of compact type equipped with 
the $G$-invariant metric induced from the Killing form of the Lie algebra of 
$G$.  
%%$B_{\mathfrak g}\vert_{\mathfrak p\times\mathfrak p}$.  
Let $H$ be a symmetric subgroup of $G$ (i.e., 
$({\rm Fix}\,\tau)_0\subset H\subset{\rm Fix}\,\tau$ 
for some involution $\tau$ of $G$), where 
${\rm Fix}\,\tau$ is the fixed point group of $\tau$ and 
$({\rm Fix}\,\tau)_0$ is the identity component of ${\rm Fix}\,\tau$.  
The natural action of $H$ on $N$ is called a {\it Hermann action} 
(see [HPTT], [Kol]).  
Let $\theta$ be an involution of $G$ with $({\rm Fix}\,\theta)_0\subset K
\subset{\rm Fix}\,\theta$.  According to [Co], when $G$ is simple, 
we may assume that 
$\theta\circ\tau=\tau\circ\theta$ by replacing $H$ to a suitable conjugate 
group of $H$ if necessary except for the following three Hermann action:

\vspace{0.2truecm}

(i) $Sp(p+q)\curvearrowright SU(2p+2q)/S(U(2p-1)\times U(2q+1))
\quad(p\geq q+2)$, 

(ii) $U(p+q+1)\curvearrowright Spin(2p+2q+2)/Spin(2p+1)\times_{{\bf Z}_2}
Spin(2q+1)\quad(p\geq q+1)$, 

(iii) $Spin(3)\times_{{\bf Z}_2}Spin(5)\curvearrowright 
Spin(8)/\mu(Spin(3)\times_{{\bf Z}_2}Spin(5))$, 

\vspace{0.2truecm}

\noindent
where $\mu$ is the triality automorphism of $Spin(8)$.  
Here we note that we remove transitive Hermann actions.  

\vspace{0.5truecm}

\noindent
{\bf Assumption.} 
In the sequel, we assume that 
$\theta\circ\tau=\tau\circ\theta$.  Then the Hermann action $H\curvearrowright G/K$ 
is said to be {\it commutative}.  

\vspace{0.5truecm}

Let $\mathfrak g,\mathfrak k$ and $\mathfrak h$ be the Lie algebras of $G,K$ 
and $H$, respectively.  Denote the involutions of $\mathfrak g$ induced form 
$\theta$ and $\tau$ by the same symbols $\theta$ and $\tau$, respectively.  
Set $\mathfrak p:={\rm Ker}(\theta+{\rm id})$ and $\mathfrak q:=
{\rm Ker}(\tau+{\rm id})$.  The vector space $\mathfrak p$ is identified with 
$T_{eK}(G/K)$, where $e$ is the identity element of $G$.  
Denote by $B_{\mathfrak g}$ the Killing form of $\mathfrak g$.  
Give $G/K$ the $G$-invariant metric arising from 
$B_{\mathfrak g}\vert_{\mathfrak p\times\mathfrak p}$.  
Take a maximal abelian subspace $\mathfrak b$ of 
$\mathfrak p\cap\mathfrak q$.  
For each $\beta\in\mathfrak b^{\ast}$, we set 
$\mathfrak p_{\beta}:=\{X\in\mathfrak p\,\vert\,{\rm ad}(b)^2(X)=-\beta(b)^2X
\,\,(\forall\,b\in\mathfrak b)\}$ and 
$\triangle':=\{\beta\in\mathfrak b^{\ast}\setminus\{0\}\,\vert\,
\mathfrak p_{\beta}\not=\{0\}\}$.  This set $\triangle'$ is a root system.  
Let $\Pi'=\{\beta_1,\cdots,\beta_r\}$ be the simple root system of 
the positive root system $\triangle'_+$ of $\triangle'$ under a lexicographic 
ordering of $\mathfrak b^{\ast}$.  Set 
${\triangle'}^V_+:=\{\beta\in\triangle'_+\,\vert\,\mathfrak p_{\beta}\cap
\mathfrak q\not=\{0\}\}$ and ${\triangle'}^H_+:=\{\beta\in\triangle'_+\,\vert\,
\mathfrak p_{\beta}\cap\mathfrak h\not=\{0\}\}$.  
Define a subset $\widetilde C$ of $\mathfrak b$ by 
$$\begin{array}{l}
\displaystyle{\widetilde C:=\{b\in\mathfrak b\,\vert\,0<\beta(b)<\pi\,
(\forall\,\beta\in{\triangle'}^V_+),\,\,-\frac{\pi}{2}<\beta(b)<
\frac{\pi}{2}\,(\forall\,\beta\in{\triangle'}^H_+)\}.}
\end{array}$$
The closure $\overline{\widetilde C}$ of $\widetilde C$ is a simplicial 
complex.  
Set $C:={\rm Exp}(\widetilde C)$, where ${\rm Exp}$ is the exponential map 
of $G/K$ at $eK$.  
Each principal $H$-orbit passes through only one point of $C$ and 
each singular $H$-orbit passes through only one point of 
${\rm Exp}(\partial\widetilde C)$.  For each simplex $\sigma$ of 
$\overline{\widetilde C}$, only one minimal $H$-orbit through 
${\rm Exp}(\sigma)$ exists.  See proofs of Theorems A and B in [K2] 
(also [I]) about this fact.  For $\beta\in\triangle'_+$, we set 
$\beta=\sum\limits_{i=1}^rn^{\beta}_i\beta_i,\,\,
m_{\beta}:={\rm dim}\,\mathfrak p_{\beta},\,\,m^V_{\beta}:=
{\rm dim}(\mathfrak p_{\beta}\cap\mathfrak q)$ and 
$m^H_{\beta}:={\rm dim}(\mathfrak p_{\beta}\cap\mathfrak h)$.  
Let $Z_0$ be a point of $\mathfrak b$.  
We consider the following two conditions for 
$Z_0$:
$$({\rm I})\quad\left\{
\begin{array}{l}
\displaystyle{
\vspace{0.2truecm}\beta(Z_0)\equiv 0,\,\frac{\pi}{6},\,\frac{\pi}{3},\,
\frac{\pi}{2},\,\frac{2\pi}{3},\,\frac{5\pi}{6}\,\,
({\rm mod}\,\pi)\,\,\,\,(\forall\,\beta\in\triangle'_+)\,\,\,\&}\\
\hspace{0.5truecm}\displaystyle{
\sum_{\beta\in{\triangle'}^V_+\,\,{\rm s.t.}\,\,
\beta(Z_0)\equiv\,\frac{\pi}{6}\,\,({\rm mod}\,\pi)}
3n^{\beta}_im^V_{\beta}
+\sum_{\beta\in{\triangle'}^V_+\,\,{\rm s.t.}\,\,
\beta(Z_0)\equiv\,\frac{\pi}{3}\,\,({\rm mod}\,\pi)}
n^{\beta}_im^V_{\beta}}\\
\hspace{0.5truecm}\displaystyle{
+\sum_{\beta\in{\triangle'}^H_+\,\,{\rm s.t.}\,\,
\beta(Z_0)\equiv\frac{2\pi}{3}\,({\rm mod}\,\pi)}3n^{\beta}_im^H_{\beta}
+\sum_{\beta\in{\triangle'}^H_+\,\,{\rm s.t.}\,\,
\beta(Z_0)\equiv\frac{5\pi}{6}\,({\rm mod}\,\pi)}n^{\beta}_im^H_{\beta}}\\
\displaystyle{
=\sum_{\beta\in{\triangle'}^V_+\,\,{\rm s.t.}\,\,
\beta(Z_0)\equiv\,\frac{2\pi}{3}\,\,({\rm mod}\,\pi)}
n^{\beta}_im^V_{\beta}
+\sum_{\beta\in{\triangle'}^V_+\,\,{\rm s.t.}\,\,
\beta(Z_0)\equiv\,\frac{5\pi}{6}\,\,({\rm mod}\,\pi)}
3n^{\beta}_im^V_{\beta}}\\
\hspace{0.5truecm}\displaystyle{
+\sum_{\beta\in{\triangle'}^H_+\,\,{\rm s.t.}\,\,
\beta(Z_0)\equiv\frac{\pi}{6}\,({\rm mod}\,\pi)}n^{\beta}_im^H_{\beta}
+\sum_{\beta\in{\triangle'}^H_+\,\,{\rm s.t.}\,\,
\beta(Z_0)\equiv\frac{\pi}{3}\,({\rm mod}\,\pi)}3n^{\beta}_im^H_{\beta}}\\
\hspace{8.8truecm}\displaystyle{(i=1,\cdots,r).}
\end{array}\right.$$
and
$$({\rm II})\quad\left\{
\begin{array}{l}
\displaystyle{
\beta(Z_0)\equiv 0,\,\frac{\pi}{4},\,\frac{\pi}{2},\,\frac{3\pi}{4}\,\,
({\rm mod}\,\pi)\,\,\,\,(\forall\,\beta\in\triangle'_+)\,\,\,\&}\\
\hspace{0.5truecm}\displaystyle{
\sum_{\beta\in{\triangle'}^V_+\,\,{\rm s.t.}\,\,
\beta(Z_0)\equiv\frac{\pi}{4}\,({\rm mod}\,\pi)}n^{\beta}_im^V_{\beta}
+\sum_{\beta\in{\triangle'}^H_+\,\,{\rm s.t.}\,\,
\beta(Z_0)\equiv\frac{3\pi}{4}\,({\rm mod}\,\pi)}n^{\beta}_im^H_{\beta}}\\
\displaystyle{
=\sum_{\beta\in{\triangle'}^V_+\,\,{\rm s.t.}\,\,
\beta(Z_0)\equiv\frac{3\pi}{4}\,({\rm mod}\,\pi)}n^{\beta}_im^V_{\beta}
+\sum_{\beta\in{\triangle'}^H_+\,\,{\rm s.t.}\,\,
\beta(Z_0)\equiv\frac{\pi}{4}\,({\rm mod}\,\pi)}n^{\beta}_im^H_{\beta}}\\
\hspace{8.8truecm}\displaystyle{(i=1,\cdots,r).}
\end{array}\right.$$
Denote by $L$ the isotropy group of $H$ at ${\rm Exp}\,Z_0$.  
%For simplicity, we set $L:=H_{Z_0}$.  
Denote by $\mathfrak h$ (resp. $\mathfrak l$) the Lie algebra of $H$ 
(resp. $L$) and $B_{\mathfrak g}$ the Killing form of $\mathfrak g$.  
Also, denote by $g_I$ the induced metric on the submanifold $M$ in $G/K$ 
%and $\rho^S_{Z_0}$ the slice representation of the orbit $M$ at 
%${\rm Exp}\,Z_0$.  Under the identification of the normal space 
%$T^{\perp}_{{\rm Exp}\,Z_0}(M)$ of $M$ at ${\rm Exp}\,Z_0$ and 
%$\mathfrak m^{\perp}:=
%{\rm Ad}(\exp(-Z_0))
%((\exp\,Z_0)_{\ast}^{-1}%(T^{\perp}_{{\rm Exp}\,Z_0}(M)))(\subset\mathfrak p\cap
%\mathfrak q)$, $\rho^S_{Z_0}$ is identified with $(\rho_{G/H}\vert_L)
%\vert_{\mathfrak m^{\perp}}$, where 
%$(\rho_{G/H}\vert_L)\vert_{\mathfrak m^{\perp}}$ means the representaion 
%of $L$ with representation space $\mathfrak m^{\perp}$ defined by 
%$((\rho_{G/H}\vert_L)\vert_{\mathfrak m^{\perp}})({\it l}))(w)
%:=(\rho_{G/H}({\it l}))(w)\,\,\,({\it l}\in L,\,w\in\mathfrak m^{\perp})$ 
%in terms of the isotropy representation $\rho_{G/H}$ of $G/H$.  
and $\nabla^{\perp}$ the normal connection of the submanifold $M$.  
In the case where $(\mathfrak h,\mathfrak l)$ admits a reductive decomposition 
$\mathfrak h=\mathfrak l+\mathfrak m$, we denote the canonical connection of 
the principal $L$-bundle $\pi:H\to H/L(=M)$ with respect to this reductive 
decomposition by $\omega_{\mathfrak m}$.  
Let $F^{\perp}(M)$ be the normal frame bundle of $M$.  Define a map 
$\eta:H\to F^{\perp}(M)$ by $\eta(h)=h_{\ast}u_0$ ($h\in H$), 
where $u_0$ is an arbitrary fixed element of $F^{\perp}(M)_{{\rm Exp}\,Z_0}$, 
where $F^{\perp}(M)_{{\rm Exp}\,Z_0}$ is the fibre of $F^{\perp}(M)$ over 
${\rm Exp}\,Z_0$.  This map $\eta$ is an embedding.  
By identifying $H$ with $\eta(H)$, we regard $\pi:H\to H/L(=M)$ as a 
subbundle of $F^{\perp}(M)$.  
Denote by the same symbol $\omega_{\mathfrak m}$ the connection of 
$F^{\perp}(M)$ induced from $\omega_{\mathfrak m}$ and 
$\nabla^{\omega_{\mathfrak m}}$ the linear connection on $T^{\perp}M$ 
associated with $\omega_{\mathfrak m}$.  

In this paper, we prove the following results for the orbit 
$M=H({\rm Exp}\,Z_0)$ of the Hermann action $H\curvearrowright G/K$.  

\vspace{0.5truecm}

\noindent
{\bf Theorem A.} {\sl If $Z_0$ satisfies the condition (I) or (II), then 
the orbit $M$ is a minimal submanifold satisfying 
the following conditions:

$({\rm i})\,\,$ $(\mathfrak h,\mathfrak l)$ admits a reductive decomposition 
$\mathfrak h=\mathfrak l+\mathfrak m$ such that 
$B_{\mathfrak g}(\mathfrak l,\mathfrak m)=0$,

$({\rm ii})\,\,$ $\nabla^{\perp}=\nabla^{\omega_{\mathfrak m}}$ holds.  

\noindent
Also, $\displaystyle{\mathop{\cap}_{v\in T^{\perp}_xM}{\rm Ker}\,A_v}$ 
is equal to 
$$\begin{array}{r}
\displaystyle{g_{0\ast}(\mathfrak z_{\mathfrak p\cap\mathfrak h}(\mathfrak b))
+\sum_{\beta\in{\triangle'}^V_+\,\,{\rm s.t.}\,\,\beta(Z_0)
\equiv\frac{\pi}{2}\,\,({\rm mod}\,\pi)}g_{0\ast}
(\mathfrak p_{\beta}\cap\mathfrak q)}\\
\displaystyle{+\sum_{\beta\in{\triangle'}^H_+\,\,{\rm s.t.}\,\,\beta(Z_0)
\equiv 0\,\,({\rm mod}\,\pi)}g_{0\ast}(\mathfrak p_{\beta}\cap\mathfrak h),}
\end{array}$$
where $\mathfrak z_{\mathfrak p\cap\mathfrak h}(\mathfrak b)$ is the centralizer of 
$\mathfrak b$ in $\mathfrak p\cap\mathfrak h$.}

\vspace{0.5truecm}

Let $M$ be a submanifold in a Riemannian manifold $N$. If, for any unit normal 
vector $v$, the spectrum of the shape operator $A_v$ is invariant with respect 
to the $(-1)$-multiple (with considering the multiplicities), then $M$ is 
called an {\it austere submanifold}.  
By using Theorem A, we can show the following fact.  

\vspace{0.5truecm}

\noindent
{\bf Theorem B.} {\sl Assume that $Z_0$ satisfies the condition (I) or (II).  
If 
%${\triangle'}^V_+={\triangle'}^H_+(=\triangle'_+)$ and furthermore if 
$m^V_{\beta}=m^H_{\beta}$ for all $\beta\in\triangle'_+$ and if 
$Z_0$ satisfies 
$\beta(Z_0)\equiv 0,\,\frac{\pi}{4},\,\frac{\pi}{2},\,\frac{3\pi}{4}\,\,
({\rm mod}\,\pi)$ for all $\beta\in\triangle'_+$, then 
the orbit $M$ is an austere submanifold 
satisfying the conditions {\rm(i)} and {\rm (ii)} in Theorem A.}

\vspace{0.5truecm}

\noindent
{\it Remark 1.1.} The austere orbits of the commutative Hermann actions were 
classified in [I].  

\vspace{0.5truecm}

Also, we can show the following facts.  

\vspace{0.5truecm}

\noindent
{\bf Theorem C.} {\sl Assume that $Z_0$ satisfies the condition (I).  
In particular, if ${\triangle'}^V_+\cap{\triangle'}^H_+=\emptyset$, 
if $\beta(Z_0)\equiv 0,\,\frac{\pi}{3},\,\frac{2\pi}{3}\,\,
({\rm mod}\,\pi)$ for all $\beta\in{\triangle'}^V_+$ 
and if $\beta(Z_0)\equiv\frac{\pi}{6},\,\frac{\pi}{2},\,\frac{5\pi}{6}\,\,
({\rm mod}\,\pi)$ for all $\beta\in{\triangle'}^H_+$, then $M$ is 
a minimal submanifold satisfying 
the conditions {\rm(i)}, {\rm (ii)} in Theorem A.  
Furthermore, if the cohomogeneity of the $H$-action is equal to the rank of 
$G/K$, then $(g_I)_{eL}=\frac34B_{\mathfrak g}\vert_{\mathfrak m\times\mathfrak m}$ and 
$\displaystyle{\mathop{\cap}_{v\in T^{\perp}_xM}{\rm Ker}\,A_v
=\{0\}}$ hold.}

\vspace{0.5truecm}

\noindent
{\bf Theorem D.} {\sl Assume that $Z_0$ satisfies the condition (I).  
In particular, if ${\triangle'}^V_+\cap{\triangle'}^H_+=\emptyset$, 
if $\beta(Z_0)\equiv 0,\,\frac{\pi}{6},\,\frac{5\pi}{6}\,\,
({\rm mod}\,\pi)$ for all $\beta\in{\triangle'}^V_+$ 
and if $\beta(Z_0)\equiv\frac{\pi}{3},\,\frac{\pi}{2},\,\frac{2\pi}{3}\,\,
({\rm mod}\,\pi)$ for all $\beta\in{\triangle'}^H_+$, then 
$M$ is 
a minimal submanifold satisfying 
the conditions {\rm(i)}, {\rm (ii)} in Theorem A.  
Furthermore, if the cohomogeneity of the $H$-action is equal to the rank of 
$G/K$, then $(g_I)_{eL}=\frac14B_{\mathfrak g}\vert_{\mathfrak m\times\mathfrak m}$ and 
$\displaystyle{\mathop{\cap}_{v\in T^{\perp}_xM}{\rm Ker}\,A_v=\{0\}}$ hold.}

\vspace{0.5truecm}

\noindent
{\bf Theorem E.} {\sl Assume that $Z_0$ satisfies the condition (II).  
In particular, if ${\triangle'}^V_+\cap{\triangle'}^H_+=\emptyset$, 
if $\beta(Z_0)\equiv 0,\,\frac{\pi}{4},\,\frac{3\pi}{4}\,\,
({\rm mod}\,\pi)$ for all $\beta\in{\triangle'}^V_+$ 
and if $\beta(Z_0)\equiv\frac{\pi}{4},\,\frac{\pi}{2},\,\frac{3\pi}{4}\,\,
({\rm mod}\,\pi)$ for all $\beta\in{\triangle'}^H_+$, then 
$M$ is a minimal submanifold satisfying 
the conditions {\rm(i)}, {\rm (ii)} in Theorem A.  
Furthermore, if the cohomogeneity of the $H$-action is equal to the rank of 
$G/K$, then $(g_I)_{eL}=\frac12B_{\mathfrak g}\vert_{\mathfrak m\times\mathfrak m}$ and 
$\displaystyle{\mathop{\cap}_{v\in T^{\perp}_xM}{\rm Ker}\,A_v=\{0\}}$ hold.}

\vspace{0.5truecm}

\noindent
{\bf Theorem F.} {\sl If ${\triangle'}^V_+\cap{\triangle'}^H_+=\emptyset$, 
if $\beta(Z_0)\equiv 0,\,\frac{\pi}{2}\,\,({\rm mod}\,\pi)$ for all 
$\beta\in{\triangle'}_+$, then $M$ is a totally geodesic submanifold 
satisfying the conditions {\rm(i)}, {\rm (ii)} in Theorem A.  
Furthermore, if the cohomogeneity of the $H$-action is equal to the rank of 
$G/K$, then $(g_I)_{eL}=B_{\mathfrak g}\vert_{\mathfrak m\times\mathfrak m}$ holds.}

\vspace{0.5truecm}

\noindent
{\it Remark 1.2.} (i) If $H=K$ then we have ${\triangle'}^H_+=\emptyset$ and hence 
${\triangle'}^V_+\cap{\triangle'}^H_+=\emptyset$.  

(ii) In Theorems C$\sim$F, when $G$ is simple, 
there exists an inner automorphism $\rho$ of $G$ with $\rho(K)=H$ 
by Proposition 4.39 of [I].  
%Also, according to the proof of the proposition, 
%%when the cohomogeneity of the $H$-action is equal to the rank of $G/K$ 
%%(even if $G$ is not simple), 
%it is shown that, when $(G,K)=(G'\times G',\triangle G')$ and $G'$ is simple, 
%so is the $H$-action.  
%After all, when $G/K$ is irreducible, so is the $H$-action.  
%However, when $G/K$ is reducible, so is not necessarily the $H$-action (see Example ?).  

\vspace{0.5truecm}

In the final section, we give examples of Hermann actions 
$H\curvearrowright G/K$ and $Z_0\in\mathfrak b$ as in Theorems B, C and F.  

\section{Basic notions and facts} 
In this section, we recall some basic notions and facts.  

\vspace{0.3truecm}

\noindent
{\bf Shape operators of orbits of Hermann actions}

\vspace{0.1truecm}

Let $H\curvearrowright G/K$ be a Hermann action and $\theta$ (resp. $\tau$) 
an involution of $G$ with $({\rm Fix}\,\theta)_0\subset K\subset{\rm Fix}\,
\theta$ (resp. $({\rm Fix}\,\tau)_0\subset H\subset{\rm Fix}\,\tau$).  Assume 
that $\theta\circ\tau=\tau\circ\theta$.  Let $\mathfrak k,\,\mathfrak p,\,
\mathfrak h,\,\mathfrak q,\,\mathfrak b,\,\mathfrak p_{\beta},\,\triangle',\,
{\triangle'}^V_+$ and ${\triangle'}^H_+$ be as in Introduction.  
Fix $Z_0\in\mathfrak b$.  
Set $M:=H({\rm Exp}\,Z_0)$ and $g_0:=\exp\,Z_0$, where 
${\rm Exp}$ is the exponential map of $G/K$ at $eK$ and 
$\exp$ is the exponential map of $G$.  
Set 
$${\triangle'}^V_{Z_0}:=\{\beta\in{\triangle'}^V_+\,\vert\,
\beta(Z_0)\equiv 0\,\,({\rm mod}\,\pi)\}$$
and 
$${\triangle'}^H_{Z_0}:=\{\beta\in{\triangle'}^H_+\,\vert\,
\beta(Z_0)\equiv\frac{\pi}{2}\,\,({\rm mod}\,\pi)\}.$$
Denote by $A$ the shape tensor of $M$.  
The tangent space $T_{{\rm Exp}\,Z_0}M$ of $M$ at ${\rm Exp}\,Z_0$ is given by 
$$T_{{\rm Exp}\,Z_0}M=g_{0\ast}\left(
\mathfrak z_{\mathfrak p\cap\mathfrak h}(\mathfrak b)+
\sum_{\beta\in{\triangle'}^V_+\setminus{\triangle'}^V_{Z_0}}
(\mathfrak p_{\beta}\cap\mathfrak q)
+\sum_{\beta\in{\triangle'}^H_+\setminus
{\triangle'}^H_{Z_0}}(\mathfrak p_{\beta}\cap\mathfrak h)\right)\leqno{(2.1)}$$
and hence 
$$T^{\perp}_{{\rm Exp}\,Z_0}M=g_{0\ast}\left(\mathfrak b
+\sum_{\beta\in{\triangle'}^V_{Z_0}}(\mathfrak p_{\beta}\cap\mathfrak q)
+\sum_{\beta\in{\triangle'}^H_{Z_0}}(\mathfrak p_{\beta}\cap\mathfrak h)
\right).\leqno{(2.2)}$$
%where ${\triangle'}^V_{Z_0}$ and ${\triangle'}^H_{Z_0}$ are as stated in 
%Introduction.  
Denote by $L$ the isotropy group of the $H$-action at ${\rm Exp}\,Z_0$.  
The slice representation $\rho^S_{Z_0}:L\to GL(T^{\perp}_{{\rm Exp}\,Z_0}M)$ of the 
$H$-action at ${\rm Exp}\,Z_0$ is given by $\rho^S_{Z_0}(h)=h_{\ast{\rm Exp}\,Z_0}
\vert_{T^{\perp}_{{\rm Exp}\,Z_0}M}\,\,\,\,(h\in H_{Z_0})$.  Then we have 
$\displaystyle{\mathop{\cup}_{h\in H_{Z_0}}\rho^S_{Z_0}(h)
(g_{0\ast}\mathfrak b)=T^{\perp}_{{\rm Exp}\,Z_0}M}$ and 
$$\begin{array}{l}
\displaystyle{A_{\rho_{Z_0}^S(h)(g_{0\ast}v)}
\vert_{\rho^S_{Z_0}(h)(g_{0\ast}
(\mathfrak z_{\mathfrak p\cap\mathfrak h}(\mathfrak b)))}=0,}\\
\displaystyle{A_{\rho_{Z_0}^S(h)(g_{0\ast}v)}\vert_{\rho^S_{Z_0}(h)
(g_{0\ast}(\mathfrak p_{\beta}\cap\mathfrak q))}
=-\frac{\beta(v)}{\tan\beta(Z_0)}{\rm id}\,\,\,
(\beta\in{\triangle'}^V_+\setminus{\triangle'}^V_{Z_0}),}\\
\displaystyle{A_{\rho_{Z_0}^S(h)(g_{0\ast}v)}\vert_{\rho^S_{Z_0}(h)
(g_{0\ast}(\mathfrak p_{\beta}\cap\mathfrak h))}
=\beta(v)\tan\beta(Z_0){\rm id}\,\,\,
(\beta\in{\triangle'}^H_+\setminus{\triangle'}^H_{Z_0}),}
\end{array}\leqno{(2.3)}$$
where $h\in L$ and $v\in\mathfrak b$.  

\vspace{0.2truecm}

\noindent
{\bf The canonical connection}

\vspace{0.1truecm}

Let $H/L$ be a reductive homogeneous space and $\mathfrak h=\mathfrak l
+\mathfrak m$ be a reductive decomposition (i.e., 
$[\mathfrak l,\mathfrak m]\subset\mathfrak m$), where $\mathfrak h$ 
(resp. $\mathfrak l$) is the Lie algebra of $H$ (resp. $L$).  
Also, let $\pi:P\to H/L$ be a principal $G$-bundle, where $G$ is a Lie group.  
Assume that $H$ acts on $P$ as $\pi(h\cdot u)=h\cdot\pi(u)$ for any 
$u\in P$ and any $h\in H$.  
Then there uniquely exists a connection $\omega$ of 
$P$ such that, for any $X\in\mathfrak m$ and any $u\in P$, 
$t\mapsto(\exp\,tX)(u)$ is a horizontal curve with respect to 
$\omega$, where $\exp$ is the exponential map of $H$.  This connection 
$\omega$ is called the {\it canonical connection} of $P$ associated with 
the reductive decomposition $\mathfrak h=\mathfrak l+\mathfrak m$.  

\section{Proof of Theorems A$\sim$F}
In this section, we shall first prove Theorems A$\sim$F.  
We use the notations in Introduction.  
Let $H\curvearrowright G/K$ be a Hermann action and $Z_0$ be an element of 
$\mathfrak b$.  Set $M:=H({\rm Exp}\,Z_0)$.  

\vspace{0.3truecm}

\noindent
{\it Proof of Theorem A.} 
Denote by ${\cal H}$ the mean curvature vector of $M$.  
From $(2.1)$ and $(2.3)$, we have 
$$\begin{array}{l}
\displaystyle{\langle{\cal H}_{{\rm Exp}\,Z_0},\rho^S_{Z_0}(h)(g_{0\ast}v)
\rangle
=-\sum_{i=1}^r\sum_{\beta\in{\triangle'}^V_+\setminus
{\triangle'}^V_{Z_0}}\frac{n^{\beta}_im^V_{\beta}}
{\tan\beta(Z_0)}\beta_i(v)}\\
\hspace{4.15truecm}
\displaystyle{+\sum_{i=1}^r\sum_{\beta\in{\triangle'}^H_+\setminus
{\triangle'}^H_{Z_0}}n^{\beta}_im^H_{\beta}\tan\beta(Z_0)\beta_i(v)}
\end{array}$$
for any $v\in\mathfrak b$ and any $h\in L$.  
Hence, ${\cal H}_{{\rm Exp}\,Z_0}$ vanishes if and only if 
the following relations hold:
$$\begin{array}{r}
\displaystyle{\sum_{\beta\in{\triangle'}^V_+\setminus
{\triangle'}^V_{Z_0}}\frac{n^{\beta}_im^V_{\beta}}{\tan\beta(Z_0)}
=\sum_{\beta\in{\triangle'}^H_+\setminus{\triangle'}^H_{Z_0}}
n^{\beta}_im^H_{\beta}\tan\beta(Z_0)}\\
(i=1,\cdots,r).
\end{array}\leqno{(3.1)}$$
Since $Z_0$ satisfies the condition (I) or (II) in Theorem A, 
$(3.1)$ holds, that is, ${\cal H}_{{\rm Exp}\,Z_0}$ vanishes.  
Therefore $M$ is minimal.  

Next we shall show that there exists a reductive 
decomposition $\mathfrak h=\mathfrak l+\mathfrak m$ with 
$B_{\mathfrak g}(\mathfrak l,\mathfrak m)=0$.  
Easily we have 
$$\mathfrak l=\mathfrak z_{\mathfrak k\cap\mathfrak h}(\mathfrak b)
+\sum_{\beta\in{\triangle'}^V_{Z_0}}
(\mathfrak k_{\beta}\cap\mathfrak h)
+\sum_{\beta\in{\triangle'}^H_{Z_0}}(\mathfrak p_{\beta}\cap\mathfrak h).
\leqno{(3.2)}$$
Define a subspace $\mathfrak m$ of $\mathfrak h$ by 
$$\mathfrak m:=\mathfrak z_{\mathfrak p\cap\mathfrak h}(\mathfrak b)
+\sum_{\beta\in{\triangle'}^V_+\setminus{\triangle'}^V_{Z_0}}
(\mathfrak k_{\beta}\cap\mathfrak h)
+\sum_{\beta\in{\triangle'}^H_+\setminus{\triangle'}^H_{Z_0}}
(\mathfrak p_{\beta}\cap\mathfrak h).\leqno{(3.3)}$$
Easily we can show that $\mathfrak h=\mathfrak l+\mathfrak m$ is a reductive 
decomposition and that $B_{\mathfrak g}(\mathfrak l,\mathfrak m)=0$.  

Next we shall show that $\nabla^{\omega_{\mathfrak m}}=\nabla^{\perp}$.  
Take $v\in\mathfrak b\,(\subset g_{0\ast}^{-1}T^{\perp}_{{\rm Exp}\,Z_0}M)$.  
Set $g_s:=\exp(1-s)Z_0$.  Let $Z:[0,1]\to \mathfrak b$ be a $C^{\infty}$-curve 
such that $Z(0)=Z_0$ and that $Z((0,1])$ is contained in a fundamental domain 
of the Coxeter group associated with the principal $H$-orbit at 
an intersection point of the orbit and $\mathfrak b$.  
Set $M_s:=H({\rm Exp}\,Z(1-s))$ 
($0\leq s\leq1$).  Denote by $A^s$ the shape tensor of $M_s$ and 
$\widetilde{\nabla}$ the Levi-Civita connection of $G/K$.  
Let ${\widetilde v}^s$ be the $H$-equivariant normal vector field of $M_s$ 
($0\leq s<1$) arising from $g_{s\ast}v$.  
Since $M_s$ ($0\leq s<1$) is a principal orbit of a Hermann (hence hyperpolar) 
action, ${\widetilde v}^s$ is well-defined and it is a parallel normal vector 
field with respect to $\nabla^{\perp}$.  
Take $X\in\mathfrak k_{\beta}\cap\mathfrak h\,(\subset\mathfrak m)$ 
($\beta\in{\triangle'}^V_+\setminus{\triangle'}^V_{Z_0}$).  
Then, by using $(2.3)$,  we have 
$$\widetilde{\nabla}_{X^{\ast}_{{\rm Exp}\,Z(1-s)}}{\widetilde v}^s=
-A^s_vX^{\ast}_{{\rm Exp}\,Z(1-s)}=\frac{\beta(v)}{\tan\beta(Z_0)}
X^{\ast}_{{\rm Exp}\,Z(1-s)}.$$
and hence 
$$\widetilde{\nabla}_{X^{\ast}_{{\rm Exp}\,Z_0}}
(\exp\,tX)_{\ast{\rm Exp}(Z_0)}(v)
=\lim_{s\to 1-0}\widetilde{\nabla}_{X^{\ast}_{{\rm Exp}\,Z(1-s)}}
{\widetilde v}^s
=\frac{\beta(v)}{\tan\beta(Z_0)}X^{\ast}_{{\rm Exp}\,Z_0}
\in T_{{\rm Exp}\,Z_0}M.$$
Hence we obtain $\nabla^{\perp}_{X^{\ast}_{{\rm Exp}\,Z_0}}
(\exp\,tX)_{\ast{\rm Exp}(Z_0)}(v)=0$.  
Take $Y\in\mathfrak p_{\beta}\cap\mathfrak h\,(\subset\mathfrak m)$ 
($\beta\in{\triangle'}^H_+\setminus{\triangle'}^H_{Z_0}$).  
Then, by using $(2.3)$,  we have 
$$\widetilde{\nabla}_{Y^{\ast}_{{\rm Exp}\,Z(1-s)}}{\widetilde v}^s=
-A^s_vY^{\ast}_{{\rm Exp}\,Z(1-s)}=-\beta(v)\tan\beta(Z_0)
Y^{\ast}_{{\rm Exp}\,Z(1-s)}.$$
and hence 
$$\begin{array}{l}
\displaystyle{\widetilde{\nabla}_{Y^{\ast}_{{\rm Exp}\,Z_0}}
(\exp\,tY)_{\ast{\rm Exp}\,Z_0}(v)
=\lim_{s\to 1-0}\widetilde{\nabla}_{Y^{\ast}_{{\rm Exp}\,Z(1-s)}}
{\widetilde v}^s}\\
\hspace{2.8truecm}
\displaystyle{=-\beta(v)\tan\beta(Z_0)Y^{\ast}_{{\rm Exp}\,Z_0}
\in T_{{\rm Exp}\,Z_0}M.}
\end{array}$$
Hence we obtain $\nabla^{\perp}_{Y^{\ast}_{{\rm Exp}\,Z_0}}
(\exp\,tY)_{\ast{\rm Exp}(Z_0)}(v)=0$.  
Therefore, it follows from the arbitrariness of $X,\,Y$ and $\beta$ that 
$t\mapsto(\exp\,t\hat X)_{\ast{\rm Exp}\,Z_0}(v)$ is $\nabla^{\perp}$-parallel 
along $t\mapsto(\exp\,t\hat X)({\rm Exp}\,Z_0)$ for any 
$\hat X\in\mathfrak m$.  
Take any $h\in L$.  Similarly we can show that 
$t\mapsto(\exp\,t\hat X)_{\ast{\rm Exp}\,Z_0}(\rho^S_{Z_0}(h))(g_{0\ast}v))$ 
is $\nabla^{\perp}$-parallel along 
$t\mapsto(\exp\,t\hat X)({\rm Exp}\,Z_0)$ for any $\hat X\in\mathfrak m$.  
Note that this fact has been showed in [IST] in different 
method.  On the other hand, it follows from the definition of $\omega$ that 
$t\mapsto(\exp\,t\hat X)_{\ast{\rm Exp}\,Z_0}(\rho^S_{Z_0}(h)(g_{0\ast}v))$ is 
$\nabla^{\omega_{\mathfrak m}}$-parallel along 
$t\mapsto(\exp\,t\hat X)({\rm Exp}\,Z_0)$ for any $\hat X\in\mathfrak m$.  
Therefore we obtain $\nabla^{\perp}=\nabla^{\omega_{\mathfrak m}}$.  
The statement for 
$\displaystyle{\mathop{\cap}_{v\in T^{\perp}_xM}{\rm Ker}\,A_v}$ follows from 
$(2.3)$ directly.  
\begin{flushright}q.e.d.\end{flushright}

\vspace{0.3truecm}

Next we prove Theorem B.  

\vspace{0.3truecm}

\noindent
{\it Proof of Theorem B.} This statement of this theorem follows from $(2.3)$ 
directly.  
\begin{flushright}q.e.d.\end{flushright}

\vspace{0.5truecm}

Next we prove Theorems C$\sim$F.  

\vspace{0.5truecm}

\noindent
{\it Proof of Theorems C$\sim$F.} 
Define a diffeomorphism $\psi:H/L\to M$ by $\psi(hL):=h\cdot{\rm Exp}\,Z_0$ 
($h\in H$).  Next we shall show that $(\psi^{\ast}g_I)_{eL}
=cB_{\mathfrak g}\vert_{\mathfrak m\times\mathfrak m}$, where 
$$c=\left\{
\begin{array}{ll}
\displaystyle{\frac34} & ({\rm in}\,\,{\rm case}\,\,{\rm of}\,\,
{\rm Theorems}\,\,{\rm C})\\
\displaystyle{\frac14} & ({\rm in}\,\,{\rm case}\,\,{\rm of}\,\,
{\rm Theorem}\,\,{\rm D})\\
\displaystyle{\frac12} & ({\rm in}\,\,{\rm case}\,\,{\rm of}\,\,
{\rm Theorem}\,\,{\rm E})\\
\displaystyle{1} & ({\rm in}\,\,{\rm case}\,\,{\rm of}\,\,
{\rm Theorem}\,\,{\rm F}).
\end{array}\right.$$
In the sequel, we omit the notation $\psi^{\ast}$.  
For each $X\in\mathfrak m(=T_{eL}(H/L)=T_{{\rm Exp}\,Z_0}M)$, 
denote by $X^{\ast}$ the Killing field on $M$ associated with $X$, that is, 
$X^{\ast}_p:=\frac{d}{dt}\vert_{t=0}
(\exp\,tX)(p)$\newline
$(p\in M)$.  From the definition of $\psi$, we have 
$\psi_{\ast eL}X=X^{\ast}_{{\rm Exp}\,Z_0}$.  
Take $S_{\beta_1}\in\mathfrak f_{\beta_1}\cap\mathfrak h$ 
($\beta_1\in{\triangle'_+}^H\setminus{\triangle'_{Z_0}}^H$) 
and $\hat S_{\beta_2}\in\mathfrak p_{\beta_2}\cap\mathfrak h$ 
($\beta_2\in{\triangle'_+}^V\setminus{\triangle'_{Z_0}}^V$).  
Let $T_{\beta_1}$ be the element of $\mathfrak p_{\beta_1}\cap\mathfrak q$ 
such that ${\rm ad}(b)(S_{\beta_1})=\beta_1(b)T_{\beta_1}$ for any 
$b\in\mathfrak b$.  Then we have 
$$\psi_{\ast eL}(S_{\beta_1})=(S_{\beta_1}^{\ast})_{{\rm Exp}\,Z_0}
=-\sin\beta_1(Z_0)(\exp\,Z_0)_{\ast}(T_{\beta_1})\leqno{(3.4)}$$
and 
$$\psi_{\ast eL}(\hat S_{\beta_2})=(\hat S_{\beta_2}^{\ast})_{{\rm Exp}\,Z_0}
=\cos\beta_2(Z_0)(\exp\,{Z_0})_{\ast}(\hat S_{\beta_2}).\leqno{(3.5)}$$
Hence, since $H$ and $Z_0$ is as in Theorems C$\sim$F, we have 
$(g_I)_{eL}(S_{\beta_1},S_{\beta_1})=c
B_{\mathfrak g}(S_{\beta_1},S_{\beta_1})$ and 
$(g_I)_{eL}(\hat S_{\beta_2},\hat S_{\beta_2})
=cB_{\mathfrak g}(\hat S_{\beta_2},\hat S_{\beta_2})$.  
If the cohomogeneity of the $H$-action is equal to the rank of $G/K$, then 
we have $\mathfrak z_{\mathfrak p\cap\mathfrak h}(\mathfrak b)=0$.  
Therefore we obtain $(g_I)_{eL}=cB_{\mathfrak g}\vert_{\mathfrak m\times\mathfrak m}$.  
Also, in Theorems C$\sim$E, 
$\displaystyle{\mathop{\cap}_{v\in T^{\perp}_xM}{\rm Ker}\,A_v}=\{0\}$ follows from 
the statement for $\displaystyle{\mathop{\cap}_{v\in T^{\perp}_xM}{\rm Ker}\,A_v}$ in 
Theorem A directly.  \hspace{6.5truecm}q.e.d.
%\begin{flushright}q.e.d.\end{flushright}

\section{Examples} 
In this section, we give examples of a Hermann action 
$H\curvearrowright G/K$ and $Z_0\in\widetilde C$ as in Theorems B, C and F.  
We use the notations in Introduction.  

\vspace{0.5truecm}

\noindent
{\it Example 1.} We consider the isotropy action of $SU(3n+3)/SO(3n+3)$.  
Then we have $\triangle_+=\triangle'_+={\triangle'}^V_+$ (which is of 
$({\mathfrak a}_{3n+2})$-type) and ${\triangle'}^H_+=\emptyset$.  
Let $\Pi=\{\beta_1,\cdots,\beta_{3n+2}\}$ be 
a simple root system of $\triangle'_+$, where we order 
$\beta_1,\cdots,\beta_{3n+2}$ as 
the Dynkin diagram of $\triangle'_+$ is as in Fig. 1, 
%%$\displaystyle{\mathop{\circ}_{\beta_1}-\mathop{\circ}_{\beta_2}-\cdots
%%-\mathop{\circ}_{\beta_{3n+2}}}$, 
$\triangle'_+=\{\beta_i+\cdots+\beta_j\,\vert\,1\leq i,j\leq 3n+2\}$.  
For any $\beta\in\triangle'_+$, we have $m_{\beta}=1$.  
Let $Z_0$ be the point of $\mathfrak b$ defined by 
$\beta_{n+1}(Z_0)=\beta_{2n+2}(Z_0)=\frac{\pi}{3}$ and $\beta_i(Z_0)=0$ 
($i\in\{1,\cdots,3n+2\}\setminus\{n+1,2n+2\}$).  
Clearly we have $m_{\beta}^V=1$, $m_{\beta}^H=0$ 
and $\beta(Z_0)\equiv 0,\frac{\pi}{3}$ or $\frac{2\pi}{3}$ 
(${\rm mod}\,\,\pi$) for any $\beta\in\triangle'_+$.  
For simplicity, set 
$\beta_{ij}:=\beta_i+\cdots+\beta_j\,\,(1\leq i\leq j\leq 3n+2)$.  
Easily we can show 
$$\begin{array}{l}
\hspace{0.5truecm}
\displaystyle{\{\beta\in{\triangle'}^V_+\,\vert\,\beta(Z_0)\equiv\frac{\pi}{3}\,\,
({\rm mod}\,\pi)\}}\\
\displaystyle{=\{\beta_{ij}\,\vert\,
1\leq i\leq n+1\leq j<2n+2,\,\,{\rm or}\,\,
n+1<i\leq 2n+2\leq j\leq 3n+2\}}
\end{array}$$
and 
$$\begin{array}{l}
\hspace{0.5truecm}
\displaystyle{\{\beta\in{\triangle'}^V_+\,\vert\,\beta(Z_0)\equiv\frac{2\pi}{3}\,\,
({\rm mod}\,\pi)\}}\\
\displaystyle{=\{\beta_{ij}\,\vert\,1\leq i\leq n+1, 2n+2\leq j\leq 3n+2\}.}
\end{array}$$
From these facts, it follows that the condition (I) holds.  
Thus $Z_0$ is as in the statement of Theorem C.  
Also, it is easy to show that $M$ is not austere.  

\vspace{0.5truecm}

\centerline{
%WinTpicVersion3.08
\unitlength 0.1in
\begin{picture}( 11.7000,  1.8000)(  3.3000, -7.2000)
% CIRCLE 2 0 3 0
% 4 440 590 480 620 600 590 730 590
% 
\special{pn 8}%
\special{ar 440 590 50 50  0.0000000 6.2831853}%
% CIRCLE 2 0 3 0
% 4 650 590 690 620 810 590 940 590
% 
\special{pn 8}%
\special{ar 650 590 50 50  0.0000000 6.2831853}%
% CIRCLE 2 0 3 0
% 4 1450 590 1490 620 1610 590 1740 590
% 
\special{pn 8}%
\special{ar 1450 590 50 50  0.0000000 6.2831853}%
% LINE 2 0 3 0
% 2 500 590 600 590
% 
\special{pn 8}%
\special{pa 500 590}%
\special{pa 600 590}%
\special{fp}%
% LINE 2 0 3 0
% 2 700 590 760 590
% 
\special{pn 8}%
\special{pa 700 590}%
\special{pa 760 590}%
\special{fp}%
% LINE 2 0 3 0
% 2 1330 590 1400 590
% 
\special{pn 8}%
\special{pa 1330 590}%
\special{pa 1400 590}%
\special{fp}%
% LINE 2 2 3 0
% 2 810 590 1240 590
% 
\special{pn 8}%
\special{pa 810 590}%
\special{pa 1240 590}%
\special{dt 0.045}%
% STR 2 0 3 0
% 3 330 600 330 700 1 0
% {\scriptsize$\beta_1$}
\put(3.3000,-7.0000){\makebox(0,0)[lt]{{\scriptsize$\beta_1$}}}%
% STR 2 0 3 0
% 3 570 600 570 700 1 0
% {\scriptsize$\beta_2$}
\put(5.7000,-7.0000){\makebox(0,0)[lt]{{\scriptsize$\beta_2$}}}%
% STR 2 0 3 0
% 3 1360 620 1360 720 1 0
% {\scriptsize$\beta_{3n+2}$}
\put(13.6000,-7.2000){\makebox(0,0)[lt]{{\scriptsize$\beta_{3n+2}$}}}%
\end{picture}%
\hspace{0truecm}}

\vspace{0.5truecm}

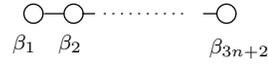
\centerline{{\bf Figure 1.}}

\vspace{0.5truecm}

\noindent
{\it Example 2.} 
We consider the isotropy action of $SU(6n+6)/Sp(3n+3)$.  
Then we have $\triangle_+=\triangle'_+={\triangle'}^V_+$ (which is of 
$({\mathfrak a}_{3n+2})$-type) and ${\triangle'}^H_+=\emptyset$.  
Let $\Pi=\{\beta_1,\cdots,\beta_{3n+2}\}$ be 
a simple root system of $\triangle'_+$, where we order 
$\beta_1,\cdots,\beta_{3n+2}$ as above.  We have $m_{\beta}=4$ for any 
$\beta\in\triangle'_+$.  
Let $Z_0$ be the point of the closure of $\mathfrak b$ defined by 
$\beta_{n+1}(Z_0)=\beta_{2n+2}(Z_0)=\frac{\pi}{3}$ and $\beta_i(Z_0)=0$ 
($i\in\{1,\cdots,3n+2\}\setminus\{n+1,2n+2\}$).  
Clearly we have $m_{\beta}^V=4$, $m_{\beta}^H=0$ 
and $\beta(Z_0)\equiv 0,\frac{\pi}{3}$ or $\frac{2\pi}{3}$ 
(${\rm mod}\,\,\pi$) for any $\beta\in\triangle'_+$.  
For simplicity, set 
$\beta_{ij}:=\beta_i+\cdots+\beta_j\,\,(1\leq i\leq j\leq 3n+2)$.  
Easily we can show 
$$\begin{array}{l}
\hspace{0.5truecm}
\displaystyle{\{\beta\in{\triangle'}^V_+\,\vert\,\beta(Z_0)\equiv\frac{\pi}{3}\,\,
({\rm mod}\,\pi)\}}\\
\displaystyle{=\{\beta_{ij}\,\vert\,
1\leq i\leq n+1\leq j<2n+2,\,\,{\rm or}\,\,
n+1<i\leq 2n+2\leq j\leq 3n+2\}}
\end{array}$$
and 
$$\begin{array}{l}
\hspace{0.5truecm}
\displaystyle{\{\beta\in{\triangle'}^V_+\,\vert\,\beta(Z_0)\equiv\frac{2\pi}{3}\,\,
({\rm mod}\,\pi)\}}\\
\displaystyle{=\{\beta_{ij}\,\vert\,1\leq i\leq n+1, 2n+2\leq j\leq 3n+2\}.}
\end{array}$$
From these facts, it follows that the condition (I) holds.  
Thus $Z_0$ is as in the statement of Theorem C.  
Also, it is easy to show that $M$ is not austere.  

\vspace{0.5truecm}

\noindent
{\it Example 3.} 
We consider the isotropy action of $SU(3)/S(U(1)\times U(2))$ 
($2$-dimensional complex projective space).  
Then we have $\triangle_+=\triangle'_+={\triangle'}^V_+=\{\beta,2\beta\}$ and 
${\triangle'}^H_+=\emptyset$, $m_{\beta}=2$ and $m_{2\beta}=1$.  
Let $Z_0$ be the point of $\mathfrak b$ defined by 
$\beta(Z_0)=\frac{\pi}{3}$.  Clearly $Z_0$ satisfies the condition (I).  
Thus $Z_0$ is as in the statement of Theorem C.  
Also, it is easy to show that $M$ is not austere.  

\vspace{0.5truecm}

\noindent
{\it Example 4.} 
We consider the isotropy action of $Sp(3n+2)/U(3n+2)$.  
Then we have $\triangle_+=\triangle'_+={\triangle'}^V_+$ (which is of 
$({\mathfrak c}_{3n+2})$-type) and ${\triangle'}^H_+=\emptyset$.  
Let $\Pi=\{\beta_1,\cdots,\beta_{3n+2}\}$ be 
a simple root system of $\triangle'_+$, where we we order 
$\beta_1,\cdots,\beta_{3n+2}$ as 
the Dynkin diagram of $\triangle'_+$ is as in Fig. 2.  
%%$\displaystyle{\mathop{\circ}_{\beta_1}-\mathop{\circ}_{\beta_2}-\cdots
%%-\mathop{\circ}_{\beta_{3n+1}}\Leftarrow\mathop{\circ}_{\beta_{3n+2}}}$.  
We have $m_{\beta}=1$ for any 
$\beta\in\triangle'_+$.  
Let $Z_0$ be the point of $\mathfrak b$ defined by 
$\beta_{n+1}(Z_0)=\beta_{2n+2}(Z_0)=\beta_{3n+2}(Z_0)=\frac{\pi}{3}$ and 
$\beta_i(Z_0)=0$ ($i\in\{1,\cdots,3n+2\}\setminus\{n+1,2n+2,3n+2\}$).  
Clearly we have $m_{\beta}^V=1$, $m_{\beta}^H=0$ 
and $\beta(Z_0)\equiv 0,\frac{\pi}{3}$ or $\frac{2\pi}{3}$ 
(${\rm mod}\,\,\pi$) for any $\beta\in\triangle'_+$.  
For simplicity, set 
$\beta_{ij}:=\beta_i+\cdots+\beta_j\,\,(1\leq i\leq j\leq 3n+2)$, 
$\widehat{\beta}_i:=2(\beta_i+\cdots+\beta_{3n+1})+\beta_{3n+2}$ 
and $\widehat{\beta}_{ij}:=\beta_i+\cdots+\beta_{j-1}
+2(\beta_j+\cdots+\beta_{3n+1})+\beta_{3n+2}\,\,(1\leq i<j\leq 3n+1)$.  
Easily we can show 
$$\begin{array}{l}
\hspace{0.5truecm}
\displaystyle{\{\beta\in{\triangle'}^V_+\,\vert\,\beta(Z_0)\equiv\frac{\pi}{3}\,\,
({\rm mod}\,\pi)\}}\\
\displaystyle{=\{\beta_{ij}\,\vert\,
1\leq i\leq n+1\leq j<2n+2\,\,{\rm or}\,\,
n+1<i\leq 2n+2\leq j<3n+2}\\
\hspace{1.25truecm}
\displaystyle{\,\,{\rm or}\,\,2n+3\leq i\leq j=3n+2\}}\\
\hspace{0.5truecm}
\displaystyle{\cup
\{\widehat{\beta}_i\,\vert\,2n+3\leq i\leq 3n+1\}}\\
\hspace{0.5truecm}
\displaystyle{\cup
\{\widehat{\beta}_{ij}\,\vert\,2n+3\leq i<j\leq 3n+1\,\,{\rm or}\,\,
1\leq i\leq n+1<j\leq 2n+2\}}
\end{array}$$
and 
$$\begin{array}{l}
\hspace{0.5truecm}
\displaystyle{\{\beta\in{\triangle'}^V_+\,\vert\,\beta(Z_0)\equiv\frac{2\pi}{3}\,\,
({\rm mod}\,\pi)\}}\\
\displaystyle{=\{\beta_{ij}\,\vert\,
"1\leq i\leq n+1\,\,\&\,\,2n+2\leq j\leq 3n+1"\,\,{\rm or}\,\,}\\
\hspace{1.25truecm}
\displaystyle{"n+2\leq i\leq 2n+2\,\,\&\,\,j=3n+2"\}}\\
\hspace{0.5truecm}
\displaystyle{\cup
\{\widehat{\beta}_i\,\vert\,1\leq i\leq n+1\}}\\
\hspace{0.5truecm}
\displaystyle{\cup
\{\widehat{\beta}_{ij}\,\vert\,1\leq i<j\leq n+1
\,\,{\rm or}\,\,n+2\leq i\leq 2n+2<j\leq 3n+1\}.}
\end{array}$$
From these facts, it follows that the condition (I) holds.  
Thus $Z_0$ is as in the statement of Theorem C.  
Also, it is easy to show that $M$ is not austere.  

\vspace{0.5truecm}

\centerline{
%WinTpicVersion3.08
\unitlength 0.1in
\begin{picture}( 14.0000,  1.6000)( 10.2000, -8.4000)
% CIRCLE 2 0 3 0
% 4 1140 730 1180 760 1300 730 1430 730
% 
\special{pn 8}%
\special{ar 1140 730 50 50  0.0000000 6.2831853}%
% CIRCLE 2 0 3 0
% 4 1350 730 1390 760 1510 730 1640 730
% 
\special{pn 8}%
\special{ar 1350 730 50 50  0.0000000 6.2831853}%
% CIRCLE 2 0 3 0
% 4 2150 730 2190 760 2310 730 2440 730
% 
\special{pn 8}%
\special{ar 2150 730 50 50  0.0000000 6.2831853}%
% LINE 2 0 3 0
% 2 1200 730 1300 730
% 
\special{pn 8}%
\special{pa 1200 730}%
\special{pa 1300 730}%
\special{fp}%
% LINE 2 0 3 0
% 2 1400 730 1460 730
% 
\special{pn 8}%
\special{pa 1400 730}%
\special{pa 1460 730}%
\special{fp}%
% LINE 2 0 3 0
% 2 2030 730 2100 730
% 
\special{pn 8}%
\special{pa 2030 730}%
\special{pa 2100 730}%
\special{fp}%
% LINE 2 2 3 0
% 2 1510 730 1940 730
% 
\special{pn 8}%
\special{pa 1510 730}%
\special{pa 1940 730}%
\special{dt 0.045}%
% CIRCLE 2 0 3 0
% 4 2370 730 2330 760 2080 730 2210 730
% 
\special{pn 8}%
\special{ar 2370 730 50 50  0.0000000 6.2831853}%
% LINE 2 0 3 0
% 2 2320 710 2230 710
% 
\special{pn 8}%
\special{pa 2320 710}%
\special{pa 2230 710}%
\special{fp}%
% LINE 2 0 3 0
% 2 2330 760 2240 760
% 
\special{pn 8}%
\special{pa 2330 760}%
\special{pa 2240 760}%
\special{fp}%
% LINE 2 0 3 0
% 2 2200 730 2270 680
% 
\special{pn 8}%
\special{pa 2200 730}%
\special{pa 2270 680}%
\special{fp}%
% LINE 2 0 3 0
% 2 2230 780 2230 780
% 
\special{pn 8}%
\special{pa 2230 780}%
\special{pa 2230 780}%
\special{fp}%
% LINE 2 0 3 0
% 2 2260 790 2190 730
% 
\special{pn 8}%
\special{pa 2260 790}%
\special{pa 2190 730}%
\special{fp}%
% STR 2 0 3 0
% 3 1020 720 1020 820 1 0
% {\scriptsize$\beta_1$}
\put(10.2000,-8.2000){\makebox(0,0)[lt]{{\scriptsize$\beta_1$}}}%
% STR 2 0 3 0
% 3 1280 730 1280 830 1 0
% {\scriptsize$\beta_2$}
\put(12.8000,-8.3000){\makebox(0,0)[lt]{{\scriptsize$\beta_2$}}}%
% STR 2 0 3 0
% 3 2310 740 2310 840 1 0
% {\scriptsize$\beta_{3n+2}$}
\put(23.1000,-8.4000){\makebox(0,0)[lt]{{\scriptsize$\beta_{3n+2}$}}}%
% STR 2 0 3 0
% 3 1910 730 1910 830 1 0
% {\scriptsize$\beta_{3n+1}$}
\put(19.1000,-8.3000){\makebox(0,0)[lt]{{\scriptsize$\beta_{3n+1}$}}}%
\end{picture}%
}

\vspace{0.5truecm}

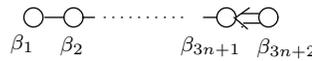
\centerline{{\bf Figure 2.}}

\vspace{0.5truecm}

%Next we shall give examples of $H$ and $Z_0$ as in Theorems B, D, E and F.  First 
By refering Tables 1 and 2 in [K2], we shall list up Hermann actions of cohomogeneity two 
on irreducible symmetric spaces of compact type and rank two satisfying 

\vspace{0.2truecm}

\centerline{
${\rm (i)}$ $m_{\beta}^V=m_{\beta}^H\,\,(\forall\,\beta\in\triangle'_+)$ 
$\qquad$ or $\qquad$ 
${\rm (ii)}$ ${\triangle'}^V_+\cap{\triangle'}^H_+=\emptyset$.  
}

\vspace{0.2truecm}

\noindent
All of such Hermann actions satisfying $({\rm i})$ are as in Table 1.  
In Table 1, $\displaystyle{\mathop{\beta}_{(m)}}$ means 
$m_{\beta}^V=m_{\beta}^H=m$.  
All of such Hermann actions satisfying $({\rm ii})$ are the dual actions (see Table 3) of 
Hermann actions on symmetric spaces of non-compact type as in Table 2.  
In Table 3, $\displaystyle{\mathop{\beta}_{(m)}}$ means 
$m_{\beta}^V$ or $m_{\beta}^H$ is equal to $m$.  
Since the Hermann actions in Table 2 are commutative, 
so are also the Hermann actions in Table 3.  
Also, since ${\triangle'}^V_+\cap{\triangle'}^H_+=\emptyset$ as in Table 3 and $G/K$ is 
irreducible, there exists an inner automorphism $\rho$ of $G$ with $\rho(K)=H$ by 
Proposition 4.39 in [I].  
According to the proof of the proposition, $\rho$ is given explicitly by 
$\rho={\rm Ad}_G(\exp\,b)$, where ${\rm Ad}_G$ is the adjoint representation of $G$ 
and $b$ is the element of $\mathfrak b$ satisfying 
$$(\beta_1(b),\beta_2(b))=\left\{
\begin{tabular}{ll}
$\displaystyle{(0,\frac{\pi}{2})}$ & (in case of (1),(2),(3),(4),(6),(9),(10),(11))\\
$\displaystyle{(\frac{\pi}{2},0)}$ & (in case of (5),(7))\\
$\displaystyle{(\frac{\pi}{2},\frac{\pi}{2})}$ & (in case of (8)).\\
\end{tabular}\right.$$

\vspace{0.3truecm}

$$\begin{tabular}{|c|c|}
\hline
{\scriptsize$H\curvearrowright G/K$} & 
{\scriptsize${\triangle'}^V_+={\triangle'}^H_+$}\\
\hline
{\scriptsize$SO(6)\curvearrowright SU(6)/Sp(3)$} & 
{\scriptsize$\displaystyle{\{\mathop{\beta_1}_{(2)},\mathop{\beta_2}_{(2)},
\mathop{\beta_1+\beta_2}_{(2)}\}}$}\\
\hline
{\scriptsize$SO(2)^2\times SO(3)^2\curvearrowright(SO(5)\times SO(5))/SO(5)$} 
& {\scriptsize$\displaystyle{\{\mathop{\beta_1}_{(1)},
\mathop{\beta_2}_{(1)},\mathop{\beta_1+\beta_2}_{(1)},
\mathop{2\beta_1+\beta_2}_{(1)}\}}$}\\
\hline
{\scriptsize$SU(2)^2\cdot SO(2)^2\curvearrowright(Sp(2)\times Sp(2))/Sp(2)$} & 
{\scriptsize$\displaystyle{\{\mathop{\beta_1}_{(1)},\mathop{\beta_2}_{(1)},
\mathop{\beta_1+\beta_2}_{(1)},\mathop{2\beta_1+\beta_2}_{(1)}\}}$}\\
\hline
{\scriptsize$Sp(4)\curvearrowright E_6/F_4$} & 
{\scriptsize$\displaystyle{\{\mathop{\beta_1}_{(4)},\mathop{\beta_2}_{(4)},
\mathop{\beta_1+\beta_2}_{(4)}\}}$}\\
\hline
{\scriptsize$SU(2)^4\curvearrowright(G_2\times G_2)/G_2$} & 
{\scriptsize$\displaystyle{\{\mathop{\beta_1}_{(1)},\mathop{\beta_2}_{(1)},
\mathop{\beta_1+\beta_2}_{(1)},\mathop{2\beta_1+\beta_2}_{(1)},
\mathop{3\beta_1+\beta_2}_{(1)},\mathop{3\beta_1+2\beta_2}_{(1)}\}}$}\\
\hline
\end{tabular}$$

\vspace{0.3truecm}

\centerline{\bf Table 1.}

\vspace{0.3truecm}

$$\begin{tabular}{|c|c|}
\hline
(1) & {\scriptsize$SO_0(1,2)\curvearrowright SL(3,{\Bbb R})/SO(3)$}\\
\hline
(2) & {\scriptsize$Sp(1,2)\curvearrowright SU^{\ast}(6)/Sp(3)$}\\
\hline
(3) & {\scriptsize$U(2,3)\curvearrowright SO^{\ast}(10)/U(5)$}\\
\hline
(4) & {\scriptsize$SO_0(2,3)\curvearrowright SO(5,{\Bbb C})/SO(5)$}\\
\hline
(5) & {\scriptsize$U(1,1)\curvearrowright Sp(2,{\Bbb R})/U(2)$}\\
\hline
(6) & {\scriptsize$Sp(2,{\Bbb R})\curvearrowright Sp(2,{\Bbb C})/Sp(2)$}\\
\hline
(7) & {\scriptsize$Sp(1,1)\curvearrowright Sp(2,{\Bbb C})/Sp(2)$}\\
\hline
(8) & {\scriptsize$SO^{\ast}(10)\cdot U(1)\curvearrowright 
E_6^{-14}/Spin(10)\cdot U(1)$}\\
\hline
(9) & {\scriptsize$F_4^{-20}\curvearrowright E_6^{-26}/F_4$}\\
\hline
(10) & {\scriptsize$SL(2,{\Bbb R})\times SL(2,{\Bbb R})
\curvearrowright G_2^2/SO(4)$}\\
\hline
(11) & {\scriptsize$G_2^2\curvearrowright G_2^{\bf C}/G_2$}\\
\hline
\end{tabular}$$

\vspace{0.3truecm}

\centerline{{\bf Table 2.}}

\vspace{0.5truecm}

$$\begin{tabular}{|c|c|c|c|}
\hline
 & {\scriptsize$H\curvearrowright G/K$} & {\scriptsize${\triangle'}^V_+$} 
& {\scriptsize${\triangle'}^H_+$}\\
\hline
(1) & {\scriptsize$SO_0(1,2)^{\ast}\curvearrowright SU(3)/SO(3)$} & 
{\scriptsize$\displaystyle{\{\mathop{\beta_1}_{(1)}\}}$} & 
{\scriptsize$\displaystyle{\{\mathop{\beta_2}_{(1)},
\mathop{\beta_1+\beta_2}_{(1)}\}}$}\\
\hline
(2) & {\scriptsize$Sp(1,2)^{\ast}\curvearrowright SU(6)/Sp(3)$} & 
{\scriptsize$\displaystyle{\{\mathop{\beta_1}_{(4)}\}}$} & 
{\scriptsize$\displaystyle{\{\mathop{\beta_2}_{(4)},
\mathop{\beta_1+\beta_2}_{(4)}\}}$}\\
\hline
(3) & {\scriptsize$U(2,3)^{\ast}\curvearrowright SO(10)/U(5)$} & 
{\scriptsize$\displaystyle{\{\mathop{\beta_1}_{(4)},\mathop{2\beta_1}_{(1)},
\mathop{2\beta_1+2\beta_2}_{(1)}\}}$} & 
{\scriptsize$\displaystyle{\{\mathop{\beta_2}_{(4)},
\mathop{\beta_1+\beta_2}_{(4)},
\mathop{2\beta_1+\beta_2}_{(4)}\}}$}\\
\hline
(4) & {\scriptsize$SO_0(2,3)^{\ast}\curvearrowright(SO(5)\times SO(5))/SO(5)$} & 
{\scriptsize$\displaystyle{\{\mathop{\beta_1}_{(2)}\}}$} & 
{\scriptsize$\displaystyle{\{\mathop{\beta_2}_{(2)},
\mathop{\beta_1+\beta_2}_{(2)},\mathop{2\beta_1+\beta_2}_{(2)}\}}$}\\
\hline
(5) & {\scriptsize$U(1,1)^{\ast}\curvearrowright Sp(2)/U(2)$} & 
{\scriptsize$\displaystyle{\{\mathop{\beta_2}_{(1)},
\mathop{2\beta_1+\beta_2}_{(1)}\}}$} & 
{\scriptsize$\displaystyle{\{\mathop{\beta_1}_{(1)},
\mathop{\beta_1+\beta_2}_{(1)}\}}$}\\
\hline
(6) & {\scriptsize$Sp(2,{\Bbb R})^{\ast}\curvearrowright (Sp(2)\times Sp(2))/Sp(2)$} & 
{\scriptsize$\displaystyle{\{\mathop{\beta_1}_{(2)}\}}$} & 
{\scriptsize$\displaystyle{\{\mathop{\beta_2}_{(2)},
\mathop{\beta_1+\beta_2}_{(2)},\mathop{2\beta_1+\beta_2}_{(2)}\}}$}\\
\hline
(7) & {\scriptsize$Sp(1,1)^{\ast}\curvearrowright(Sp(2)\times Sp(2))/Sp(2)$} & 
{\scriptsize$\displaystyle{\{\mathop{\beta_2}_{(2)},
\mathop{2\beta_1+\beta_2}_{(2)}\}}$} & 
{\scriptsize$\displaystyle{\{\mathop{\beta_1}_{(2)},
\mathop{\beta_1+\beta_2}_{(2)}\}}$}\\
\hline
(8) & {\scriptsize$(SO^{\ast}(10)\cdot U(1))^{\ast}\curvearrowright 
E_6/Spin(10)\cdot U(1)$} & 
{\scriptsize$\displaystyle{\{\mathop{\beta_1}_{(8)},\mathop{2\beta_1}_{(1)},
\mathop{2\beta_1+2\beta_2}_{(1)}\}}$} & 
{\scriptsize$\displaystyle{\{
\mathop{\beta_2}_{(6)},\mathop{\beta_1+\beta_2}_{(9)},
\mathop{2\beta_1+\beta_2}_{(5)}\}}$}\\
\hline
(9) & {\scriptsize$(F_4^{-20})^{\ast}\curvearrowright E_6/F_4$} & 
{\scriptsize$\displaystyle{\{\mathop{\beta_1}_{(8)}\}}$} 
& {\scriptsize$\displaystyle{\{\mathop{\beta_2}_{(8)},
\mathop{\beta_1+\beta_2}_{(8)}\}}$}\\
\hline
(10) & {\scriptsize$(SL(2,{\Bbb R})\times SL(2,{\Bbb R}))^{\ast}\curvearrowright G_2/SO(4)$} 
& {\scriptsize$\displaystyle{\{\mathop{\beta_1}_{(1)},
\mathop{3\beta_1+2\beta_2}_{(1)}\}}$} & 
{\scriptsize$\displaystyle{\{\mathop{\beta_2}_{(1)},
\mathop{\beta_1+\beta_2}_{(1)},
\mathop{2\beta_1+\beta_2}_{(1)},\mathop{3\beta_1+\beta_2}_{(1)}\}}$}\\
\hline
(11) & {\scriptsize$(G_2^2)^{\ast}\curvearrowright(G_2\times G_2)/G_2$} & 
{\scriptsize$\displaystyle{\{\mathop{\beta_1}_{(2)},
\mathop{3\beta_1+2\beta_2}_{(2)}\}}$} & 
{\scriptsize$\displaystyle{\{\mathop{\beta_2}_{(2)},
\mathop{\beta_1+\beta_2}_{(2)},
\mathop{2\beta_1+\beta_2}_{(2)},\mathop{3\beta_1+\beta_2}_{(2)}\}}$}\\
\hline
\end{tabular}$$

\vspace{0.3truecm}

\centerline{{\bf Table 3.}}

\vspace{0.3truecm}

According to Theorem B, we obtain the following fact.  

\vspace{0.3truecm}

\noindent
{\bf Proposition 4.1.} {\sl Let $H\curvearrowright G/K$ be a Hermann action 
in Table 1 and $Z_0$ an element of $\mathfrak b$ satisfying 
$(\beta_1(Z_0),\beta_2(Z_0))=(0,\frac{\pi}{4}),\,(\frac{\pi}{4},0)$ or 
$(\frac{\pi}{4},\frac{\pi}{4})$.  Then $M=H({\rm Exp}\,Z_0)$ is a 
(non-totally geodesic) austere submanifold.}

\vspace{0.25truecm}

\noindent
Denote by $Z_{(a,b)}$ the element $Z$ of $\mathfrak b$ satisfying 
$(\beta_1(Z),\beta_2(Z))=(a,b)$.  In the case where $\triangle'$ is 
of type ($\mathfrak a_2$), three points of $\mathfrak b$ as in 
Proposition 4.1 are as in Figure 3.  

\vspace{0.15truecm}

\centerline{
%\input{EXHEMF1.tex}
%WinTpicVersion3.08
\unitlength 0.1in
\begin{picture}( 38.7200, 19.3900)( -2.8200,-23.1000)
% VECTOR 2 0 3 0
% 2 1350 1729 3590 1729
% 
\special{pn 8}%
\special{pa 1350 1730}%
\special{pa 3590 1730}%
\special{fp}%
\special{sh 1}%
\special{pa 3590 1730}%
\special{pa 3524 1710}%
\special{pa 3538 1730}%
\special{pa 3524 1750}%
\special{pa 3590 1730}%
\special{fp}%
% VECTOR 2 0 3 0
% 4 2192 2181 2192 547 2192 547 2192 547
% 
\special{pn 8}%
\special{pa 2192 2182}%
\special{pa 2192 548}%
\special{fp}%
\special{sh 1}%
\special{pa 2192 548}%
\special{pa 2172 614}%
\special{pa 2192 600}%
\special{pa 2212 614}%
\special{pa 2192 548}%
\special{fp}%
\special{pa 2192 548}%
\special{pa 2192 548}%
\special{fp}%
% LINE 2 0 3 0
% 2 1974 2068 2868 618
% 
\special{pn 8}%
\special{pa 1974 2068}%
\special{pa 2868 618}%
\special{fp}%
% LINE 2 0 3 0
% 2 3134 2074 2240 625
% 
\special{pn 8}%
\special{pa 3134 2074}%
\special{pa 2240 626}%
\special{fp}%
% DOT 0 0 3 0
% 2 2553 1729 2553 1729
% 
\special{pn 20}%
\special{sh 1}%
\special{ar 2554 1730 10 10 0  6.28318530717959E+0000}%
\special{sh 1}%
\special{ar 2554 1730 10 10 0  6.28318530717959E+0000}%
% DOT 0 0 3 0
% 2 2360 1435 2360 1435
% 
\special{pn 20}%
\special{sh 1}%
\special{ar 2360 1436 10 10 0  6.28318530717959E+0000}%
\special{sh 1}%
\special{ar 2360 1436 10 10 0  6.28318530717959E+0000}%
% DOT 0 0 3 0
% 2 2748 1435 2748 1435
% 
\special{pn 20}%
\special{sh 1}%
\special{ar 2748 1436 10 10 0  6.28318530717959E+0000}%
\special{sh 1}%
\special{ar 2748 1436 10 10 0  6.28318530717959E+0000}%
% LINE 3 0 3 0
% 38 2221 1729 2553 1134 2259 1729 2573 1166 2297 1729 2593 1199 2336 1729 2613 1232 2373 1729 2633 1265 2410 1729 2652 1298 2449 1729 2670 1331 2487 1729 2690 1363 2524 1729 2710 1397 2562 1729 2730 1430 2600 1729 2749 1462 2639 1729 2769 1496 2675 1729 2789 1528 2713 1729 2807 1561 2752 1729 2827 1593 2790 1729 2846 1626 2827 1729 2866 1659 2865 1729 2886 1692 2903 1729 2905 1725
% 
\special{pn 4}%
\special{pa 2222 1730}%
\special{pa 2554 1134}%
\special{fp}%
\special{pa 2260 1730}%
\special{pa 2574 1166}%
\special{fp}%
\special{pa 2298 1730}%
\special{pa 2594 1200}%
\special{fp}%
\special{pa 2336 1730}%
\special{pa 2614 1232}%
\special{fp}%
\special{pa 2374 1730}%
\special{pa 2634 1266}%
\special{fp}%
\special{pa 2410 1730}%
\special{pa 2652 1298}%
\special{fp}%
\special{pa 2450 1730}%
\special{pa 2670 1332}%
\special{fp}%
\special{pa 2488 1730}%
\special{pa 2690 1364}%
\special{fp}%
\special{pa 2524 1730}%
\special{pa 2710 1398}%
\special{fp}%
\special{pa 2562 1730}%
\special{pa 2730 1430}%
\special{fp}%
\special{pa 2600 1730}%
\special{pa 2750 1462}%
\special{fp}%
\special{pa 2640 1730}%
\special{pa 2770 1496}%
\special{fp}%
\special{pa 2676 1730}%
\special{pa 2790 1528}%
\special{fp}%
\special{pa 2714 1730}%
\special{pa 2808 1562}%
\special{fp}%
\special{pa 2752 1730}%
\special{pa 2828 1594}%
\special{fp}%
\special{pa 2790 1730}%
\special{pa 2846 1626}%
\special{fp}%
\special{pa 2828 1730}%
\special{pa 2866 1660}%
\special{fp}%
\special{pa 2866 1730}%
\special{pa 2886 1692}%
\special{fp}%
\special{pa 2904 1730}%
\special{pa 2906 1726}%
\special{fp}%
% VECTOR 2 2 3 0
% 2 1822 1014 2352 1427
% 
\special{pn 8}%
\special{pa 1822 1014}%
\special{pa 2352 1428}%
\special{dt 0.045}%
\special{sh 1}%
\special{pa 2352 1428}%
\special{pa 2312 1370}%
\special{pa 2310 1394}%
\special{pa 2288 1402}%
\special{pa 2352 1428}%
\special{fp}%
% VECTOR 2 2 3 0
% 2 3378 1383 2756 1427
% 
\special{pn 8}%
\special{pa 3378 1384}%
\special{pa 2756 1428}%
\special{dt 0.045}%
\special{sh 1}%
\special{pa 2756 1428}%
\special{pa 2824 1442}%
\special{pa 2810 1424}%
\special{pa 2822 1402}%
\special{pa 2756 1428}%
\special{fp}%
% VECTOR 2 2 3 0
% 2 2293 2270 2545 1747
% 
\special{pn 8}%
\special{pa 2294 2270}%
\special{pa 2546 1748}%
\special{dt 0.045}%
\special{sh 1}%
\special{pa 2546 1748}%
\special{pa 2498 1798}%
\special{pa 2522 1796}%
\special{pa 2534 1816}%
\special{pa 2546 1748}%
\special{fp}%
% VECTOR 2 2 3 0
% 2 3390 2275 3130 1733
% 
\special{pn 8}%
\special{pa 3390 2276}%
\special{pa 3130 1734}%
\special{dt 0.045}%
\special{sh 1}%
\special{pa 3130 1734}%
\special{pa 3142 1802}%
\special{pa 3154 1782}%
\special{pa 3178 1784}%
\special{pa 3130 1734}%
\special{fp}%
% STR 2 0 3 0
% 3 2669 465 2669 541 2 0
% $(\beta_1+\beta_2)^{-1}(\frac{\pi}{2})$
\put(26.6900,-5.4100){\makebox(0,0)[lb]{$(\beta_1+\beta_2)^{-1}(\frac{\pi}{2})$}}%
% STR 2 0 3 0
% 3 3133 829 3133 906 1 0
% $(\beta_2)^{-1}(0)$
\put(31.3300,-9.0600){\makebox(0,0)[lt]{$(\beta_2)^{-1}(0)$}}%
% STR 2 0 3 0
% 3 3311 2235 3311 2310 1 0
% $(\beta_1)^{-1}(0)$
\put(33.1100,-23.1000){\makebox(0,0)[lt]{$(\beta_1)^{-1}(0)$}}%
% STR 2 0 3 0
% 3 3413 1247 3413 1322 1 0
% $Z_{(\frac{\pi}{4},\frac{\pi}{4})}$
\put(34.1300,-13.2200){\makebox(0,0)[lt]{$Z_{(\frac{\pi}{4},\frac{\pi}{4})}$}}%
% STR 2 0 3 0
% 3 2091 2233 2091 2309 1 0
% $Z_{(0,\frac{\pi}{4})}$
\put(20.9100,-23.0900){\makebox(0,0)[lt]{$Z_{(0,\frac{\pi}{4})}$}}%
% STR 2 0 3 0
% 3 1788 923 1788 999 3 0
% $Z_{(\frac{\pi}{4},0)}$
\put(17.8800,-9.9900){\makebox(0,0)[rb]{$Z_{(\frac{\pi}{4},0)}$}}%
% STR 2 0 3 0
% 3 2496 1375 2496 1450 1 0
% $\widetilde C$
\put(24.9600,-14.5000){\makebox(0,0)[lt]{$\widetilde C$}}%
% VECTOR 2 2 3 0
% 2 3094 914 2759 790
% 
\special{pn 8}%
\special{pa 3094 914}%
\special{pa 2760 790}%
\special{dt 0.045}%
\special{sh 1}%
\special{pa 2760 790}%
\special{pa 2816 832}%
\special{pa 2810 810}%
\special{pa 2828 794}%
\special{pa 2760 790}%
\special{fp}%
% VECTOR 2 2 3 0
% 2 2651 505 2305 745
% 
\special{pn 8}%
\special{pa 2652 506}%
\special{pa 2306 746}%
\special{dt 0.045}%
\special{sh 1}%
\special{pa 2306 746}%
\special{pa 2372 724}%
\special{pa 2350 716}%
\special{pa 2348 692}%
\special{pa 2306 746}%
\special{fp}%
\end{picture}%
\hspace{5truecm}}

\vspace{0.3truecm}

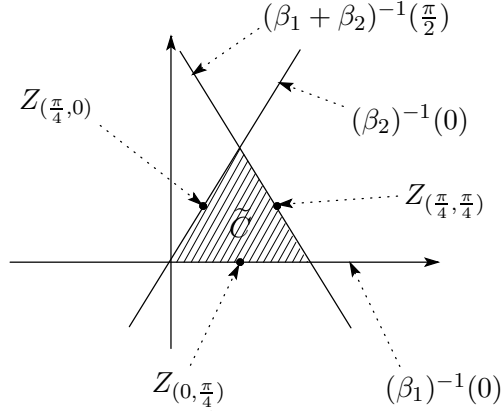
\centerline{\bf Figure 3.}

\vspace{0.25truecm}

\noindent
{\bf Proposition 4.2.} {\sl Let $H\curvearrowright G/K$ be a Hermann action 
in Table 3 and $Z_0$ an element of the closure of 
$\widetilde C(\subset\mathfrak b)$ such that $H({\rm Exp}\,Z_0)$ is minimal.  
Then, as in Tables $4\sim 13$, $Z_0$ satisfies the condition in 
Theorem C or F, or it does not satisfy the conditions in Theorems C$\sim$F.}

\vspace{0.25truecm}

\noindent
{\it Remark 4.1.} There exist exactly seven elements $Z_0$ of the closure of 
$\widetilde C(\subset\mathfrak b)$ such that $H({\rm Exp}\,Z_0)$ is minimal.  

%(ii) Hermann actions $H\curvearrowright G/K$ in Table 2 
%%other than $\rho_8(Spin(10)\cdot U(1))\curvearrowright E_6/Spin(10)\cdot U(1)$ 
%are orbit equivalent to the isotropy action of $G/K$.  

\vspace{0.15truecm}

$$\begin{tabular}{|c|c|c|c|}
\hline
{\scriptsize$(a,b)$} & {\scriptsize $Z_{(a,b)}$} & 
{\scriptsize$M=SO_0(1,2)^{\ast}({\rm Exp}\,Z_{(a,b)})$} & 
{\scriptsize ${\rm dim}\,M$}\\
\hline
{\scriptsize $(0,-\frac{\pi}{2})$} & 
{\scriptsize as in Theorem F} & {\scriptsize one-point set} & $0$\\
\hline
{\scriptsize $(0,\frac{\pi}{2})$} & 
{\scriptsize as in Theorem F} & {\scriptsize one-point set} & $0$\\
\hline
{\scriptsize $(\pi,-\frac{\pi}{2})$} & 
{\scriptsize as in Theorem F} & {\scriptsize one-point set} & $0$\\
\hline
{\scriptsize $(0,0)$} & 
{\scriptsize as in Theorem F} & {\scriptsize totally geodesic} & $2$\\
\hline
{\scriptsize $(\frac{\pi}{2},0)$} & 
{\scriptsize as in Theorem F} & {\scriptsize totally geodesic} & $2$\\
\hline
{\scriptsize $(\frac{\pi}{2},-\frac{\pi}{2})$} & 
{\scriptsize as in Theorem F} & {\scriptsize totally geodesic} & $2$\\
\hline
{\scriptsize $(\frac{\pi}{3},-\frac{\pi}{6})$} & 
{\scriptsize as in Theorem C} & {\scriptsize not austere} & $3$\\
\hline
\end{tabular}$$

\vspace{0.15truecm}

\centerline{$SO_0(1,2)^{\ast}\curvearrowright SU(3)/SO(3)$}

\centerline{$({\rm dim}\,SU(3)/SO(3)=5)$}

\vspace{0.3truecm}

\centerline{{\bf Table 4.}}

\vspace{0.5truecm}

The positions of $Z_0$'s in Table 4 are as in Figure 4.  

\vspace{0.5truecm}

\centerline{
%\input{EXHEMF2.tex}
%WinTpicVersion3.08
\unitlength 0.1in
\begin{picture}( 48.8200, 21.5800)(-12.4000,-25.4800)
% VECTOR 2 0 3 0
% 2 1372 1896 3642 1896
% 
\special{pn 8}%
\special{pa 1372 1896}%
\special{pa 3642 1896}%
\special{fp}%
\special{sh 1}%
\special{pa 3642 1896}%
\special{pa 3576 1876}%
\special{pa 3590 1896}%
\special{pa 3576 1916}%
\special{pa 3642 1896}%
\special{fp}%
% VECTOR 2 0 3 0
% 4 2590 2480 2590 643 2590 643 2590 643
% 
\special{pn 8}%
\special{pa 2590 2480}%
\special{pa 2590 644}%
\special{fp}%
\special{sh 1}%
\special{pa 2590 644}%
\special{pa 2570 710}%
\special{pa 2590 696}%
\special{pa 2610 710}%
\special{pa 2590 644}%
\special{fp}%
\special{pa 2590 644}%
\special{pa 2590 644}%
\special{fp}%
% LINE 2 0 3 0
% 2 2004 2277 2911 647
% 
\special{pn 8}%
\special{pa 2004 2278}%
\special{pa 2912 648}%
\special{fp}%
% LINE 2 0 3 0
% 2 3181 2285 2274 655
% 
\special{pn 8}%
\special{pa 3182 2286}%
\special{pa 2274 656}%
\special{fp}%
% DOT 0 0 3 0
% 2 2592 1896 2592 1896
% 
\special{pn 20}%
\special{sh 1}%
\special{ar 2592 1896 10 10 0  6.28318530717959E+0000}%
\special{sh 1}%
\special{ar 2592 1896 10 10 0  6.28318530717959E+0000}%
% DOT 0 0 3 0
% 2 2396 1566 2396 1566
% 
\special{pn 20}%
\special{sh 1}%
\special{ar 2396 1566 10 10 0  6.28318530717959E+0000}%
\special{sh 1}%
\special{ar 2396 1566 10 10 0  6.28318530717959E+0000}%
% DOT 0 0 3 0
% 2 2789 1566 2789 1566
% 
\special{pn 20}%
\special{sh 1}%
\special{ar 2790 1566 10 10 0  6.28318530717959E+0000}%
\special{sh 1}%
\special{ar 2790 1566 10 10 0  6.28318530717959E+0000}%
% LINE 3 0 3 0
% 38 2255 1896 2592 1227 2294 1896 2612 1264 2332 1896 2632 1301 2371 1896 2652 1338 2409 1896 2672 1375 2447 1896 2691 1412 2486 1896 2711 1449 2524 1896 2731 1485 2563 1896 2751 1523 2601 1896 2771 1560 2639 1896 2790 1596 2678 1896 2810 1634 2716 1896 2830 1670 2754 1896 2849 1707 2793 1896 2870 1744 2831 1896 2889 1781 2870 1896 2909 1818 2908 1896 2929 1855 2946 1896 2948 1892
% 
\special{pn 4}%
\special{pa 2256 1896}%
\special{pa 2592 1228}%
\special{fp}%
\special{pa 2294 1896}%
\special{pa 2612 1264}%
\special{fp}%
\special{pa 2332 1896}%
\special{pa 2632 1302}%
\special{fp}%
\special{pa 2372 1896}%
\special{pa 2652 1338}%
\special{fp}%
\special{pa 2410 1896}%
\special{pa 2672 1376}%
\special{fp}%
\special{pa 2448 1896}%
\special{pa 2692 1412}%
\special{fp}%
\special{pa 2486 1896}%
\special{pa 2712 1450}%
\special{fp}%
\special{pa 2524 1896}%
\special{pa 2732 1486}%
\special{fp}%
\special{pa 2564 1896}%
\special{pa 2752 1524}%
\special{fp}%
\special{pa 2602 1896}%
\special{pa 2772 1560}%
\special{fp}%
\special{pa 2640 1896}%
\special{pa 2790 1596}%
\special{fp}%
\special{pa 2678 1896}%
\special{pa 2810 1634}%
\special{fp}%
\special{pa 2716 1896}%
\special{pa 2830 1670}%
\special{fp}%
\special{pa 2754 1896}%
\special{pa 2850 1708}%
\special{fp}%
\special{pa 2794 1896}%
\special{pa 2870 1744}%
\special{fp}%
\special{pa 2832 1896}%
\special{pa 2890 1782}%
\special{fp}%
\special{pa 2870 1896}%
\special{pa 2910 1818}%
\special{fp}%
\special{pa 2908 1896}%
\special{pa 2930 1856}%
\special{fp}%
\special{pa 2946 1896}%
\special{pa 2948 1892}%
\special{fp}%
% VECTOR 2 2 3 0
% 2 3424 1471 2797 1557
% 
\special{pn 8}%
\special{pa 3424 1472}%
\special{pa 2798 1558}%
\special{dt 0.045}%
\special{sh 1}%
\special{pa 2798 1558}%
\special{pa 2866 1568}%
\special{pa 2850 1550}%
\special{pa 2860 1528}%
\special{pa 2798 1558}%
\special{fp}%
% VECTOR 2 2 3 0
% 2 2328 2505 2584 1916
% 
\special{pn 8}%
\special{pa 2328 2506}%
\special{pa 2584 1916}%
\special{dt 0.045}%
\special{sh 1}%
\special{pa 2584 1916}%
\special{pa 2540 1970}%
\special{pa 2564 1966}%
\special{pa 2576 1986}%
\special{pa 2584 1916}%
\special{fp}%
% STR 2 0 3 0
% 3 2550 475 2550 560 2 0
% $(\beta_1+\beta_2)^{-1}(\frac{\pi}{2})$
\put(25.5000,-5.6000){\makebox(0,0)[lb]{$(\beta_1+\beta_2)^{-1}(\frac{\pi}{2})$}}%
% STR 2 0 3 0
% 3 3180 825 3180 910 1 0
% $(\beta_2)^{-1}(-\frac{\pi}{2})$
\put(31.8000,-9.1000){\makebox(0,0)[lt]{$(\beta_2)^{-1}(-\frac{\pi}{2})$}}%
% STR 2 0 3 0
% 3 1360 2115 1360 2200 1 0
% $(\beta_1)^{-1}(0)$
\put(13.6000,-22.0000){\makebox(0,0)[lt]{$(\beta_1)^{-1}(0)$}}%
% STR 2 0 3 0
% 3 3460 1285 3460 1370 1 0
% $Z_{(\frac{\pi}{2},0)}$
\put(34.6000,-13.7000){\makebox(0,0)[lt]{$Z_{(\frac{\pi}{2},0)}$}}%
% STR 2 0 3 0
% 3 2123 2463 2123 2548 1 0
% $Z_{(0,0)}$
\put(21.2300,-25.4800){\makebox(0,0)[lt]{$Z_{(0,0)}$}}%
% STR 2 0 3 0
% 3 2000 1435 2000 1520 3 0
% $Z_{(\frac{\pi}{2},-\frac{\pi}{2})}$
\put(20.0000,-15.2000){\makebox(0,0)[rb]{$Z_{(\frac{\pi}{2},-\frac{\pi}{2})}$}}%
% DOT 0 0 3 0
% 2 2210 1900 2210 1900
% 
\special{pn 20}%
\special{sh 1}%
\special{ar 2210 1900 10 10 0  6.28318530717959E+0000}%
\special{sh 1}%
\special{ar 2210 1900 10 10 0  6.28318530717959E+0000}%
% DOT 0 0 3 0
% 2 2590 1230 2590 1230
% 
\special{pn 20}%
\special{sh 1}%
\special{ar 2590 1230 10 10 0  6.28318530717959E+0000}%
\special{sh 1}%
\special{ar 2590 1230 10 10 0  6.28318530717959E+0000}%
% DOT 2 0 3 0
% 2 2960 1900 2960 1900
% 
\special{pn 8}%
\special{sh 1}%
\special{ar 2960 1900 10 10 0  6.28318530717959E+0000}%
\special{sh 1}%
\special{ar 2960 1900 10 10 0  6.28318530717959E+0000}%
% VECTOR 2 2 3 0
% 2 1860 1740 2204 1894
% 
\special{pn 8}%
\special{pa 1860 1740}%
\special{pa 2204 1894}%
\special{dt 0.045}%
\special{sh 1}%
\special{pa 2204 1894}%
\special{pa 2152 1850}%
\special{pa 2156 1872}%
\special{pa 2136 1886}%
\special{pa 2204 1894}%
\special{fp}%
% STR 2 0 3 0
% 3 1840 1705 1840 1790 3 0
% $Z_{(0,-\frac{\pi}{2})}$
\put(18.4000,-17.9000){\makebox(0,0)[rb]{$Z_{(0,-\frac{\pi}{2})}$}}%
% VECTOR 2 2 3 0
% 2 1530 2160 1640 1900
% 
\special{pn 8}%
\special{pa 1530 2160}%
\special{pa 1640 1900}%
\special{dt 0.045}%
\special{sh 1}%
\special{pa 1640 1900}%
\special{pa 1596 1954}%
\special{pa 1620 1950}%
\special{pa 1632 1970}%
\special{pa 1640 1900}%
\special{fp}%
% VECTOR 2 2 3 0
% 2 3440 2120 2970 1910
% 
\special{pn 8}%
\special{pa 3440 2120}%
\special{pa 2970 1910}%
\special{dt 0.045}%
\special{sh 1}%
\special{pa 2970 1910}%
\special{pa 3024 1956}%
\special{pa 3020 1932}%
\special{pa 3040 1920}%
\special{pa 2970 1910}%
\special{fp}%
% STR 2 0 3 0
% 3 3480 2015 3480 2100 1 0
% $Z_{(0,\frac{\pi}{2})}$
\put(34.8000,-21.0000){\makebox(0,0)[lt]{$Z_{(0,\frac{\pi}{2})}$}}%
% STR 2 0 3 0
% 3 2190 1125 2190 1210 3 0
% $Z_{(\pi,-\frac{\pi}{2})}$
\put(21.9000,-12.1000){\makebox(0,0)[rb]{$Z_{(\pi,-\frac{\pi}{2})}$}}%
% VECTOR 2 2 3 0
% 2 2020 1460 2390 1570
% 
\special{pn 8}%
\special{pa 2020 1460}%
\special{pa 2390 1570}%
\special{dt 0.045}%
\special{sh 1}%
\special{pa 2390 1570}%
\special{pa 2332 1532}%
\special{pa 2340 1556}%
\special{pa 2320 1570}%
\special{pa 2390 1570}%
\special{fp}%
% LINE 2 2 3 0
% 2 2210 1900 2770 1570
% 
\special{pn 8}%
\special{pa 2210 1900}%
\special{pa 2770 1570}%
\special{dt 0.045}%
% DOT 0 0 3 0
% 2 2590 1680 2590 1680
% 
\special{pn 20}%
\special{sh 1}%
\special{ar 2590 1680 10 10 0  6.28318530717959E+0000}%
\special{sh 1}%
\special{ar 2590 1680 10 10 0  6.28318530717959E+0000}%
% VECTOR 2 2 3 0
% 2 3440 1750 2600 1680
% 
\special{pn 8}%
\special{pa 3440 1750}%
\special{pa 2600 1680}%
\special{dt 0.045}%
\special{sh 1}%
\special{pa 2600 1680}%
\special{pa 2666 1706}%
\special{pa 2654 1684}%
\special{pa 2668 1666}%
\special{pa 2600 1680}%
\special{fp}%
% STR 2 0 3 0
% 3 3470 1585 3470 1670 1 0
% $Z_{(\frac{\pi}{3},-\frac{\pi}{6})}$
\put(34.7000,-16.7000){\makebox(0,0)[lt]{$Z_{(\frac{\pi}{3},-\frac{\pi}{6})}$}}%
% VECTOR 2 2 3 0
% 2 2520 520 2310 720
% 
\special{pn 8}%
\special{pa 2520 520}%
\special{pa 2310 720}%
\special{dt 0.045}%
\special{sh 1}%
\special{pa 2310 720}%
\special{pa 2372 690}%
\special{pa 2350 684}%
\special{pa 2344 660}%
\special{pa 2310 720}%
\special{fp}%
% VECTOR 2 2 3 0
% 2 3150 940 2840 800
% 
\special{pn 8}%
\special{pa 3150 940}%
\special{pa 2840 800}%
\special{dt 0.045}%
\special{sh 1}%
\special{pa 2840 800}%
\special{pa 2894 846}%
\special{pa 2890 822}%
\special{pa 2910 810}%
\special{pa 2840 800}%
\special{fp}%
% VECTOR 2 2 3 0
% 2 2220 1150 2580 1230
% 
\special{pn 8}%
\special{pa 2220 1150}%
\special{pa 2580 1230}%
\special{dt 0.045}%
\special{sh 1}%
\special{pa 2580 1230}%
\special{pa 2520 1196}%
\special{pa 2528 1218}%
\special{pa 2512 1236}%
\special{pa 2580 1230}%
\special{fp}%
\end{picture}%
\hspace{7.5truecm}}

\vspace{0.65truecm}

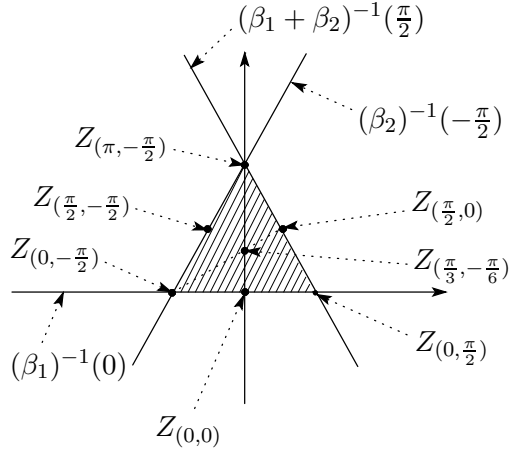
\centerline{\bf Figure 4.}

\vspace{0.5truecm}

$$\begin{tabular}{|c|c|c|c|}
\hline
{\scriptsize$(a,b)$} & {\scriptsize $Z_{(a,b)}$} & 
{\scriptsize$M=Sp(1,2)^{\ast}({\rm Exp}\,Z_{(a,b)})$} & 
{\scriptsize ${\rm dim}\,M$}\\
\hline
{\scriptsize $(0,-\frac{\pi}{2})$} & 
{\scriptsize as in Theorem F} & {\scriptsize one-point set} & $0$\\
\hline
{\scriptsize $(0,\frac{\pi}{2})$} & 
{\scriptsize as in Theorem F} & {\scriptsize one-point set} & $0$\\
\hline
{\scriptsize $(\pi,-\frac{\pi}{2})$} & 
{\scriptsize as in Theorem F} & {\scriptsize one-point set} & $0$\\
\hline
{\scriptsize $(0,0)$} & 
{\scriptsize as in Theorem F} & {\scriptsize totally geodesic} & $8$\\
\hline
{\scriptsize $(\frac{\pi}{2},0)$} & 
{\scriptsize as in Theorem F} & {\scriptsize totally geodesic} & $8$\\
\hline
{\scriptsize $(\frac{\pi}{2},-\frac{\pi}{2})$} & 
{\scriptsize as in Theorem F} & {\scriptsize totally geodesic} & $8$\\
\hline
{\scriptsize $(\frac{\pi}{3},-\frac{\pi}{6})$} & 
{\scriptsize as in Theorem C} & {\scriptsize not austere} & $12$\\
\hline
\end{tabular}$$

\vspace{0.2truecm}

\centerline{$Sp(1,2)^{\ast}\curvearrowright SU(6)/Sp(3)$}

\centerline{$({\rm dim}\,SU(6)/Sp(3)=14)$}

\vspace{0.5truecm}

\centerline{{\bf Table 5.}}

\vspace{0.5truecm}

The positions of $Z_0$'s in Table 5 are as in Figure 4.  

\vspace{0.5truecm}

$$\begin{tabular}{|c|c|c|c|}
\hline
{\scriptsize$(a,b))$} & {\scriptsize $Z_{(a,b)}$} & 
{\scriptsize$M=U(2,3)^{\ast}({\rm Exp}\,Z_{(a,b)})$} & {\scriptsize ${\rm dim}\,M$}\\
\hline
{\scriptsize $(0,\frac{\pi}{2})$} & 
{\scriptsize as in Theorem F} & {\scriptsize one-point set} & $0$\\
\hline
{\scriptsize $(0,0)$} & 
{\scriptsize as in Theorem F} & {\scriptsize totally geodesic} & $12$\\
\hline
{\scriptsize $(\frac{\pi}{2},-\frac{\pi}{2})$} & 
{\scriptsize as in Theorem F} & {\scriptsize totally geodesic} & $8$\\
\hline
{\scriptsize $(\arctan\sqrt{\frac73},\frac{\pi}{2}-\arctan\sqrt{\frac73})$} & 
{\scriptsize not as in Theorems C$\sim$F} & {\scriptsize not austere} & $14$\\
\hline
{\scriptsize $(0,\arctan\frac{1}{\sqrt{13}})$} & 
{\scriptsize not as in Theorems C$\sim$F} & {\scriptsize not austere} & $13$\\
\hline
{\scriptsize $(\arctan\frac{\sqrt5}{3},-\arctan\frac{\sqrt 5}{3})$} & 
{\scriptsize not as in Theorems C$\sim$F} & {\scriptsize not austere} & $17$\\
\hline
{\scriptsize $(a_0,b_0)$} & 
{\scriptsize not as in Theorems C$\sim$F} & {\scriptsize not austere} & $18$\\
\hline
\end{tabular}$$

\vspace{0.2truecm}

\centerline{$U(2,3)^{\ast}\curvearrowright SO(10)/U(5)$}

\centerline{$({\rm dim}\,SO(10)/U(5)=20)$}

\vspace{0.5truecm}

\centerline{{\bf Table 6.}}

\vspace{0.7truecm}

The positions of $Z_0$'s in Table 6 are as in Figure 5.  
Also, the numbers $a_0$ and $b_0$ in Table 6 are real numbers such that 
$a_0,b_0\not\equiv\frac{\pi}{6},\frac{\pi}{3},\frac{\pi}{4},\frac{3\pi}{4}\,\,
({\rm mod}\,\pi)$.  

\vspace{0.5truecm}

\centerline{
%\input{EXHEMF3.tex}
%WinTpicVersion3.08
\unitlength 0.1in
\begin{picture}( 68.3200, 17.7900)(-18.4600,-20.6400)
% VECTOR 2 0 3 0
% 2 2768 1550 4874 1550
% 
\special{pn 8}%
\special{pa 2768 1550}%
\special{pa 4874 1550}%
\special{fp}%
\special{sh 1}%
\special{pa 4874 1550}%
\special{pa 4808 1530}%
\special{pa 4822 1550}%
\special{pa 4808 1570}%
\special{pa 4874 1550}%
\special{fp}%
% VECTOR 2 0 3 0
% 4 3709 2007 3709 558 3709 558 3709 558
% 
\special{pn 8}%
\special{pa 3710 2008}%
\special{pa 3710 558}%
\special{fp}%
\special{sh 1}%
\special{pa 3710 558}%
\special{pa 3690 626}%
\special{pa 3710 612}%
\special{pa 3730 626}%
\special{pa 3710 558}%
\special{fp}%
\special{pa 3710 558}%
\special{pa 3710 558}%
\special{fp}%
% LINE 2 0 3 0
% 2 4470 1853 3305 565
% 
\special{pn 8}%
\special{pa 4470 1854}%
\special{pa 3306 566}%
\special{fp}%
% DOT 0 0 3 0
% 2 3713 1546 3713 1546
% 
\special{pn 20}%
\special{sh 1}%
\special{ar 3714 1546 10 10 0  6.28318530717959E+0000}%
\special{sh 1}%
\special{ar 3714 1546 10 10 0  6.28318530717959E+0000}%
% DOT 0 0 3 0
% 2 3905 1221 3905 1221
% 
\special{pn 20}%
\special{sh 1}%
\special{ar 3906 1222 10 10 0  6.28318530717959E+0000}%
\special{sh 1}%
\special{ar 3906 1222 10 10 0  6.28318530717959E+0000}%
% STR 2 0 3 0
% 3 3774 399 3774 455 2 0
% $(2\beta_1+\beta_2)^{-1}(\frac{\pi}{2})$
\put(37.7400,-4.5500){\makebox(0,0)[lb]{$(2\beta_1+\beta_2)^{-1}(\frac{\pi}{2})$}}%
% STR 2 0 3 0
% 3 2973 1678 2973 1744 4 0
% $(\beta_1)^{-1}(0)$
\put(29.7300,-17.4400){\makebox(0,0)[rt]{$(\beta_1)^{-1}(0)$}}%
% STR 2 0 3 0
% 3 3569 1698 3569 1764 4 0
% $Z_{(0,0)}$
\put(35.6900,-17.6400){\makebox(0,0)[rt]{$Z_{(0,0)}$}}%
% DOT 0 0 3 0
% 2 3868 1550 3868 1550
% 
\special{pn 20}%
\special{sh 1}%
\special{ar 3868 1550 10 10 0  6.28318530717959E+0000}%
\special{sh 1}%
\special{ar 3868 1550 10 10 0  6.28318530717959E+0000}%
% DOT 0 0 3 0
% 2 3709 1021 3709 1021
% 
\special{pn 20}%
\special{sh 1}%
\special{ar 3710 1022 10 10 0  6.28318530717959E+0000}%
\special{sh 1}%
\special{ar 3710 1022 10 10 0  6.28318530717959E+0000}%
% DOT 2 0 3 0
% 2 4185 1548 4185 1548
% 
\special{pn 8}%
\special{sh 1}%
\special{ar 4186 1548 10 10 0  6.28318530717959E+0000}%
\special{sh 1}%
\special{ar 4186 1548 10 10 0  6.28318530717959E+0000}%
% STR 2 0 3 0
% 3 4641 1610 4641 1667 1 0
% $Z_{(0,\frac{\pi}{2})}$
\put(46.4100,-16.6700){\makebox(0,0)[lt]{$Z_{(0,\frac{\pi}{2})}$}}%
% STR 2 0 3 0
% 3 4026 678 4026 726 2 0
% $Z_{(\frac{\pi}{2},-\frac{\pi}{2})}$
\put(40.2600,-7.2600){\makebox(0,0)[lb]{$Z_{(\frac{\pi}{2},-\frac{\pi}{2})}$}}%
% DOT 0 0 3 0
% 2 3709 1377 3709 1377
% 
\special{pn 20}%
\special{sh 1}%
\special{ar 3710 1378 10 10 0  6.28318530717959E+0000}%
\special{sh 1}%
\special{ar 3710 1378 10 10 0  6.28318530717959E+0000}%
% VECTOR 2 2 3 0
% 2 4930 1405 3851 1352
% 
\special{pn 8}%
\special{pa 4930 1406}%
\special{pa 3852 1352}%
\special{dt 0.045}%
\special{sh 1}%
\special{pa 3852 1352}%
\special{pa 3918 1376}%
\special{pa 3904 1356}%
\special{pa 3920 1336}%
\special{pa 3852 1352}%
\special{fp}%
% STR 2 0 3 0
% 3 4278 938 4278 998 1 0
% $Z_{(\arctan\sqrt{\frac 73},\frac{\pi}{2}-2\arctan\sqrt{\frac 73})}$
\put(42.7800,-9.9800){\makebox(0,0)[lt]{$Z_{(\arctan\sqrt{\frac 73},\frac{\pi}{2}-2\arctan\sqrt{\frac 73})}$}}%
% VECTOR 2 2 3 0
% 2 3504 1744 3691 1560
% 
\special{pn 8}%
\special{pa 3504 1744}%
\special{pa 3692 1560}%
\special{dt 0.045}%
\special{sh 1}%
\special{pa 3692 1560}%
\special{pa 3630 1594}%
\special{pa 3654 1598}%
\special{pa 3658 1622}%
\special{pa 3692 1560}%
\special{fp}%
% VECTOR 2 2 3 0
% 2 4007 726 3719 1007
% 
\special{pn 8}%
\special{pa 4008 726}%
\special{pa 3720 1008}%
\special{dt 0.045}%
\special{sh 1}%
\special{pa 3720 1008}%
\special{pa 3782 976}%
\special{pa 3758 970}%
\special{pa 3754 946}%
\special{pa 3720 1008}%
\special{fp}%
% VECTOR 2 2 3 0
% 2 3961 1880 3868 1579
% 
\special{pn 8}%
\special{pa 3962 1880}%
\special{pa 3868 1580}%
\special{dt 0.045}%
\special{sh 1}%
\special{pa 3868 1580}%
\special{pa 3870 1650}%
\special{pa 3884 1630}%
\special{pa 3908 1638}%
\special{pa 3868 1580}%
\special{fp}%
% STR 2 0 3 0
% 3 3793 1870 3793 1929 1 0
% $Z_{(0,\arctan\frac{1}{\sqrt{13}})}$
\put(37.9300,-19.2900){\makebox(0,0)[lt]{$Z_{(0,\arctan\frac{1}{\sqrt{13}})}$}}%
% DOT 0 0 3 0
% 2 3858 1356 3858 1356
% 
\special{pn 20}%
\special{sh 1}%
\special{ar 3858 1356 10 10 0  6.28318530717959E+0000}%
\special{sh 1}%
\special{ar 3858 1356 10 10 0  6.28318530717959E+0000}%
% STR 2 0 3 0
% 3 3374 1122 3374 1182 4 0
% $Z_{(\arctan\frac{\sqrt 5}{3},-\arctan\frac{\sqrt 5}{3})}$
\put(33.7400,-11.8200){\makebox(0,0)[rt]{$Z_{(\arctan\frac{\sqrt 5}{3},-\arctan\frac{\sqrt 5}{3})}$}}%
% STR 2 0 3 0
% 3 3309 1969 3309 2035 4 0
% $(\beta_1+\beta_2)^{-1}(0)$
\put(33.0900,-20.3500){\makebox(0,0)[rt]{$(\beta_1+\beta_2)^{-1}(0)$}}%
% VECTOR 2 2 3 0
% 2 3737 416 3364 610
% 
\special{pn 8}%
\special{pa 3738 416}%
\special{pa 3364 610}%
\special{dt 0.045}%
\special{sh 1}%
\special{pa 3364 610}%
\special{pa 3432 598}%
\special{pa 3412 586}%
\special{pa 3414 562}%
\special{pa 3364 610}%
\special{fp}%
% STR 2 0 3 0
% 3 4986 1260 4986 1356 1 0
% $Z_{(a_0,b_0)}$
\put(49.8600,-13.5600){\makebox(0,0)[lt]{$Z_{(a_0,b_0)}$}}%
% LINE 3 0 3 0
% 36 3709 1046 3723 1032 3709 1104 3750 1062 3709 1163 3778 1091 3709 1221 3805 1121 3709 1279 3832 1151 3709 1337 3860 1180 3709 1395 3887 1210 3709 1453 3915 1239 3709 1512 3942 1269 3728 1550 3970 1298 3784 1550 3997 1328 3840 1550 4025 1357 3896 1550 4052 1388 3951 1550 4080 1417 4007 1550 4107 1447 4063 1550 4134 1477 4119 1550 4162 1506 4175 1550 4189 1536
% 
\special{pn 4}%
\special{pa 3710 1046}%
\special{pa 3724 1032}%
\special{fp}%
\special{pa 3710 1104}%
\special{pa 3750 1062}%
\special{fp}%
\special{pa 3710 1164}%
\special{pa 3778 1092}%
\special{fp}%
\special{pa 3710 1222}%
\special{pa 3806 1122}%
\special{fp}%
\special{pa 3710 1280}%
\special{pa 3832 1152}%
\special{fp}%
\special{pa 3710 1338}%
\special{pa 3860 1180}%
\special{fp}%
\special{pa 3710 1396}%
\special{pa 3888 1210}%
\special{fp}%
\special{pa 3710 1454}%
\special{pa 3916 1240}%
\special{fp}%
\special{pa 3710 1512}%
\special{pa 3942 1270}%
\special{fp}%
\special{pa 3728 1550}%
\special{pa 3970 1298}%
\special{fp}%
\special{pa 3784 1550}%
\special{pa 3998 1328}%
\special{fp}%
\special{pa 3840 1550}%
\special{pa 4026 1358}%
\special{fp}%
\special{pa 3896 1550}%
\special{pa 4052 1388}%
\special{fp}%
\special{pa 3952 1550}%
\special{pa 4080 1418}%
\special{fp}%
\special{pa 4008 1550}%
\special{pa 4108 1448}%
\special{fp}%
\special{pa 4064 1550}%
\special{pa 4134 1478}%
\special{fp}%
\special{pa 4120 1550}%
\special{pa 4162 1506}%
\special{fp}%
\special{pa 4176 1550}%
\special{pa 4190 1536}%
\special{fp}%
% VECTOR 2 2 3 0
% 2 3402 1298 3700 1376
% 
\special{pn 8}%
\special{pa 3402 1298}%
\special{pa 3700 1376}%
\special{dt 0.045}%
\special{sh 1}%
\special{pa 3700 1376}%
\special{pa 3642 1340}%
\special{pa 3648 1362}%
\special{pa 3630 1378}%
\special{pa 3700 1376}%
\special{fp}%
% VECTOR 2 2 3 0
% 2 4613 1706 4203 1560
% 
\special{pn 8}%
\special{pa 4614 1706}%
\special{pa 4204 1560}%
\special{dt 0.045}%
\special{sh 1}%
\special{pa 4204 1560}%
\special{pa 4260 1602}%
\special{pa 4254 1578}%
\special{pa 4274 1564}%
\special{pa 4204 1560}%
\special{fp}%
% VECTOR 2 2 3 0
% 2 4259 1114 3905 1221
% 
\special{pn 8}%
\special{pa 4260 1114}%
\special{pa 3906 1222}%
\special{dt 0.045}%
\special{sh 1}%
\special{pa 3906 1222}%
\special{pa 3976 1222}%
\special{pa 3956 1206}%
\special{pa 3964 1184}%
\special{pa 3906 1222}%
\special{fp}%
% VECTOR 2 2 3 0
% 2 3346 2064 3709 1919
% 
\special{pn 8}%
\special{pa 3346 2064}%
\special{pa 3710 1920}%
\special{dt 0.045}%
\special{sh 1}%
\special{pa 3710 1920}%
\special{pa 3640 1926}%
\special{pa 3660 1940}%
\special{pa 3656 1962}%
\special{pa 3710 1920}%
\special{fp}%
% VECTOR 2 2 3 0
% 2 2777 1715 2871 1550
% 
\special{pn 8}%
\special{pa 2778 1716}%
\special{pa 2872 1550}%
\special{dt 0.045}%
\special{sh 1}%
\special{pa 2872 1550}%
\special{pa 2822 1598}%
\special{pa 2846 1596}%
\special{pa 2856 1618}%
\special{pa 2872 1550}%
\special{fp}%
\end{picture}%
\hspace{12truecm}}

\vspace{0.5truecm}

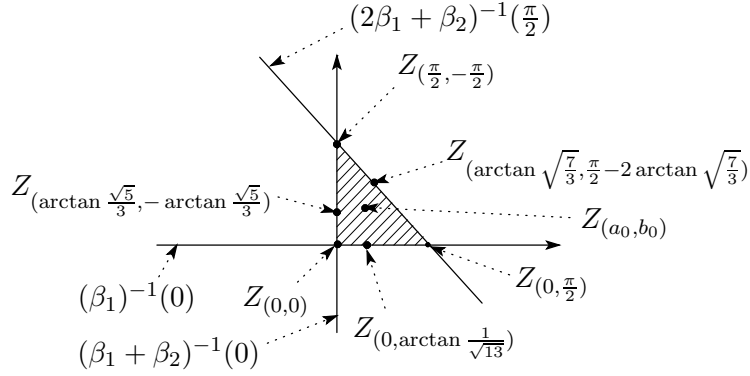
\centerline{\bf Figure 5.}

%\newpage

\vspace{0.15truecm}

$$\begin{tabular}{|c|c|c|c|}
\hline
{\scriptsize$(a,b)$} & {\scriptsize $Z_{(a,b)}$} & 
{\scriptsize$M=SO_0(2,3)^{\ast}({\rm Exp}\,Z_{(a,b)})$} & {\scriptsize ${\rm dim}\,M$}\\
\hline
{\scriptsize $(0,-\frac{\pi}{2})$} & 
{\scriptsize as in Theorem F} & {\scriptsize one-point set} & $0$\\
\hline
{\scriptsize $(0,\frac{\pi}{2})$} & 
{\scriptsize as in Theorem F} & {\scriptsize one-point set} & $0$\\
\hline
{\scriptsize $(\frac{\pi}{2},-\frac{\pi}{2})$} & 
{\scriptsize as in Theorem F} & {\scriptsize totally geodesic} & $4$\\
\hline
{\scriptsize $(0,0)$} & 
{\scriptsize as in Theorem F} & {\scriptsize totally geodesic} & $6$\\
\hline
{\scriptsize $(\arctan\sqrt 3,-\frac{\pi}{2})$} & 
{\scriptsize not as in Theorems C$\sim$F} & {\scriptsize not austere} & $6$\\
\hline
{\scriptsize $(\arctan\sqrt 3,\frac{\pi}{2}-2\arctan\sqrt 3)$} & 
{\scriptsize not as in Theorems C$\sim$F} & {\scriptsize not austere} & $6$\\
\hline
{\scriptsize $(\arctan\frac{1}{\sqrt 2},-\arctan\frac{1}{\sqrt2})$} & 
{\scriptsize not as in Theorems C$\sim$F} & {\scriptsize not austere} & $8$\\
\hline
\end{tabular}$$

\vspace{0.1truecm}

\centerline{$SO_0(2,3)^{\ast}\curvearrowright(SO(5)\times SO(5))/SO(5)$}

\centerline{$({\rm dim}\,(SO(5)\times SO(5))/SO(5)=10)$}

\vspace{0.25truecm}

\centerline{{\bf Table 7.}}

\vspace{0.3truecm}

The positions of $Z_0$'s in Table 7 are as in Figure 6.  

\vspace{0.5truecm}

\centerline{
%\input{EXHEMF4.tex}
%WinTpicVersion3.08
\unitlength 0.1in
\begin{picture}( 74.4000, 19.7100)(-21.0000,-20.6100)
% VECTOR 2 0 3 0
% 2 3080 1590 5340 1590
% 
\special{pn 8}%
\special{pa 3080 1590}%
\special{pa 5340 1590}%
\special{fp}%
\special{sh 1}%
\special{pa 5340 1590}%
\special{pa 5274 1570}%
\special{pa 5288 1590}%
\special{pa 5274 1610}%
\special{pa 5340 1590}%
\special{fp}%
% VECTOR 2 0 3 0
% 4 4090 2061 4090 566 4090 566 4090 566
% 
\special{pn 8}%
\special{pa 4090 2062}%
\special{pa 4090 566}%
\special{fp}%
\special{sh 1}%
\special{pa 4090 566}%
\special{pa 4070 634}%
\special{pa 4090 620}%
\special{pa 4110 634}%
\special{pa 4090 566}%
\special{fp}%
\special{pa 4090 566}%
\special{pa 4090 566}%
\special{fp}%
% LINE 2 0 3 0
% 2 4906 1902 3656 574
% 
\special{pn 8}%
\special{pa 4906 1902}%
\special{pa 3656 574}%
\special{fp}%
% DOT 0 0 3 0
% 2 4094 1585 4094 1585
% 
\special{pn 20}%
\special{sh 1}%
\special{ar 4094 1586 10 10 0  6.28318530717959E+0000}%
\special{sh 1}%
\special{ar 4094 1586 10 10 0  6.28318530717959E+0000}%
% DOT 0 0 3 0
% 2 4300 1250 4300 1250
% 
\special{pn 20}%
\special{sh 1}%
\special{ar 4300 1250 10 10 0  6.28318530717959E+0000}%
\special{sh 1}%
\special{ar 4300 1250 10 10 0  6.28318530717959E+0000}%
% STR 2 0 3 0
% 3 3950 202 3950 260 2 0
% $(2\beta_1+\beta_2)^{-1}(\frac{\pi}{2})$
\put(39.5000,-2.6000){\makebox(0,0)[lb]{$(2\beta_1+\beta_2)^{-1}(\frac{\pi}{2})$}}%
% STR 2 0 3 0
% 3 4070 450 4070 500 2 0
% $(\beta_2)^{-1}(-\frac{\pi}{2})$
\put(40.7000,-5.0000){\makebox(0,0)[lb]{$(\beta_2)^{-1}(-\frac{\pi}{2})$}}%
% STR 2 0 3 0
% 3 4020 1722 4020 1790 4 0
% $Z_{(0,0)}$
\put(40.2000,-17.9000){\makebox(0,0)[rt]{$Z_{(0,0)}$}}%
% DOT 0 0 3 0
% 2 4090 1044 4090 1044
% 
\special{pn 20}%
\special{sh 1}%
\special{ar 4090 1044 10 10 0  6.28318530717959E+0000}%
\special{sh 1}%
\special{ar 4090 1044 10 10 0  6.28318530717959E+0000}%
% DOT 2 0 3 0
% 2 4601 1588 4601 1588
% 
\special{pn 8}%
\special{sh 1}%
\special{ar 4602 1588 10 10 0  6.28318530717959E+0000}%
\special{sh 1}%
\special{ar 4602 1588 10 10 0  6.28318530717959E+0000}%
% VECTOR 2 2 3 0
% 2 5263 1769 4613 1596
% 
\special{pn 8}%
\special{pa 5264 1770}%
\special{pa 4614 1596}%
\special{dt 0.045}%
\special{sh 1}%
\special{pa 4614 1596}%
\special{pa 4672 1632}%
\special{pa 4666 1610}%
\special{pa 4684 1594}%
\special{pa 4614 1596}%
\special{fp}%
% STR 2 0 3 0
% 3 5290 1681 5290 1740 1 0
% $Z_{(0,\frac{\pi}{2})}$
\put(52.9000,-17.4000){\makebox(0,0)[lt]{$Z_{(0,\frac{\pi}{2})}$}}%
% STR 2 0 3 0
% 3 4550 860 4550 910 2 0
% $Z_{(\frac{\pi}{2},-\frac{\pi}{2})}$
\put(45.5000,-9.1000){\makebox(0,0)[lb]{$Z_{(\frac{\pi}{2},-\frac{\pi}{2})}$}}%
% DOT 0 0 3 0
% 2 4090 1411 4090 1411
% 
\special{pn 20}%
\special{sh 1}%
\special{ar 4090 1412 10 10 0  6.28318530717959E+0000}%
\special{sh 1}%
\special{ar 4090 1412 10 10 0  6.28318530717959E+0000}%
% STR 2 0 3 0
% 3 4750 1048 4750 1110 1 0
% $Z_{(\arctan\sqrt 3,\,\frac{\pi}{2}-2\arctan\sqrt 3)}$
\put(47.5000,-11.1000){\makebox(0,0)[lt]{$Z_{(\arctan\sqrt 3,\,\frac{\pi}{2}-2\arctan\sqrt 3)}$}}%
% STR 2 0 3 0
% 3 3300 1228 3300 1290 4 0
% $Z_{(\arctan\frac{1}{\sqrt 2},\,-\arctan\frac{1}{\sqrt 2})}$
\put(33.0000,-12.9000){\makebox(0,0)[rt]{$Z_{(\arctan\frac{1}{\sqrt 2},\,-\arctan\frac{1}{\sqrt 2})}$}}%
% LINE 2 0 3 0
% 2 3290 1890 4540 562
% 
\special{pn 8}%
\special{pa 3290 1890}%
\special{pa 4540 562}%
\special{fp}%
% VECTOR 2 2 3 0
% 2 4720 1200 4300 1250
% 
\special{pn 8}%
\special{pa 4720 1200}%
\special{pa 4300 1250}%
\special{dt 0.045}%
\special{sh 1}%
\special{pa 4300 1250}%
\special{pa 4370 1262}%
\special{pa 4354 1244}%
\special{pa 4364 1222}%
\special{pa 4300 1250}%
\special{fp}%
% DOT 0 0 3 0
% 2 3900 1240 3900 1240
% 
\special{pn 20}%
\special{sh 1}%
\special{ar 3900 1240 10 10 0  6.28318530717959E+0000}%
\special{sh 1}%
\special{ar 3900 1240 10 10 0  6.28318530717959E+0000}%
% VECTOR 2 2 3 0
% 2 3310 1370 4080 1400
% 
\special{pn 8}%
\special{pa 3310 1370}%
\special{pa 4080 1400}%
\special{dt 0.045}%
\special{sh 1}%
\special{pa 4080 1400}%
\special{pa 4014 1378}%
\special{pa 4028 1398}%
\special{pa 4014 1418}%
\special{pa 4080 1400}%
\special{fp}%
% STR 2 0 3 0
% 3 3480 888 3480 950 4 0
% $Z_{(\arctan\sqrt 3,\,-\frac{\pi}{2})}$
\put(34.8000,-9.5000){\makebox(0,0)[rt]{$Z_{(\arctan\sqrt 3,\,-\frac{\pi}{2})}$}}%
% DOT 0 0 3 0
% 2 3570 1590 3570 1590
% 
\special{pn 20}%
\special{sh 1}%
\special{ar 3570 1590 10 10 0  6.28318530717959E+0000}%
\special{sh 1}%
\special{ar 3570 1590 10 10 0  6.28318530717959E+0000}%
% STR 2 0 3 0
% 3 2810 1741 2810 1800 4 0
% $Z_{(0,-\frac{\pi}{2})}$
\put(28.1000,-18.0000){\makebox(0,0)[rt]{$Z_{(0,-\frac{\pi}{2})}$}}%
% VECTOR 2 2 3 0
% 2 2840 1860 3550 1600
% 
\special{pn 8}%
\special{pa 2840 1860}%
\special{pa 3550 1600}%
\special{dt 0.045}%
\special{sh 1}%
\special{pa 3550 1600}%
\special{pa 3482 1604}%
\special{pa 3500 1618}%
\special{pa 3494 1642}%
\special{pa 3550 1600}%
\special{fp}%
% LINE 3 0 3 0
% 34 3682 1478 4109 1051 3630 1590 4138 1082 3691 1589 4167 1113 3752 1588 4196 1144 3812 1588 4225 1175 3873 1587 4254 1206 3933 1587 4282 1238 3994 1586 4311 1269 4055 1585 4340 1300 4115 1585 4369 1331 4176 1584 4398 1362 4236 1584 4427 1393 4297 1583 4456 1424 4358 1582 4485 1455 4418 1582 4513 1487 4479 1581 4542 1518 4539 1581 4571 1549
% 
\special{pn 4}%
\special{pa 3682 1478}%
\special{pa 4110 1052}%
\special{fp}%
\special{pa 3630 1590}%
\special{pa 4138 1082}%
\special{fp}%
\special{pa 3692 1590}%
\special{pa 4168 1114}%
\special{fp}%
\special{pa 3752 1588}%
\special{pa 4196 1144}%
\special{fp}%
\special{pa 3812 1588}%
\special{pa 4226 1176}%
\special{fp}%
\special{pa 3874 1588}%
\special{pa 4254 1206}%
\special{fp}%
\special{pa 3934 1588}%
\special{pa 4282 1238}%
\special{fp}%
\special{pa 3994 1586}%
\special{pa 4312 1270}%
\special{fp}%
\special{pa 4056 1586}%
\special{pa 4340 1300}%
\special{fp}%
\special{pa 4116 1586}%
\special{pa 4370 1332}%
\special{fp}%
\special{pa 4176 1584}%
\special{pa 4398 1362}%
\special{fp}%
\special{pa 4236 1584}%
\special{pa 4428 1394}%
\special{fp}%
\special{pa 4298 1584}%
\special{pa 4456 1424}%
\special{fp}%
\special{pa 4358 1582}%
\special{pa 4486 1456}%
\special{fp}%
\special{pa 4418 1582}%
\special{pa 4514 1488}%
\special{fp}%
\special{pa 4480 1582}%
\special{pa 4542 1518}%
\special{fp}%
\special{pa 4540 1582}%
\special{pa 4572 1550}%
\special{fp}%
% VECTOR 2 2 3 0
% 2 4350 510 4430 680
% 
\special{pn 8}%
\special{pa 4350 510}%
\special{pa 4430 680}%
\special{dt 0.045}%
\special{sh 1}%
\special{pa 4430 680}%
\special{pa 4420 612}%
\special{pa 4408 632}%
\special{pa 4384 628}%
\special{pa 4430 680}%
\special{fp}%
% VECTOR 2 2 3 0
% 2 4530 870 4090 1030
% 
\special{pn 8}%
\special{pa 4530 870}%
\special{pa 4090 1030}%
\special{dt 0.045}%
\special{sh 1}%
\special{pa 4090 1030}%
\special{pa 4160 1026}%
\special{pa 4140 1012}%
\special{pa 4146 988}%
\special{pa 4090 1030}%
\special{fp}%
% VECTOR 2 2 3 0
% 2 4020 280 3760 670
% 
\special{pn 8}%
\special{pa 4020 280}%
\special{pa 3760 670}%
\special{dt 0.045}%
\special{sh 1}%
\special{pa 3760 670}%
\special{pa 3814 626}%
\special{pa 3790 626}%
\special{pa 3780 604}%
\special{pa 3760 670}%
\special{fp}%
% VECTOR 2 2 3 0
% 4 3500 1060 3880 1230 3880 1230 3880 1230
% 
\special{pn 8}%
\special{pa 3500 1060}%
\special{pa 3880 1230}%
\special{dt 0.045}%
\special{sh 1}%
\special{pa 3880 1230}%
\special{pa 3828 1186}%
\special{pa 3832 1208}%
\special{pa 3812 1222}%
\special{pa 3880 1230}%
\special{fp}%
\special{pa 3880 1230}%
\special{pa 3880 1230}%
\special{dt 0.045}%
% VECTOR 2 2 3 0
% 2 3910 1770 4090 1590
% 
\special{pn 8}%
\special{pa 3910 1770}%
\special{pa 4090 1590}%
\special{dt 0.045}%
\special{sh 1}%
\special{pa 4090 1590}%
\special{pa 4030 1624}%
\special{pa 4052 1628}%
\special{pa 4058 1652}%
\special{pa 4090 1590}%
\special{fp}%
\end{picture}%
\hspace{12.5truecm}}

\vspace{0.5truecm}

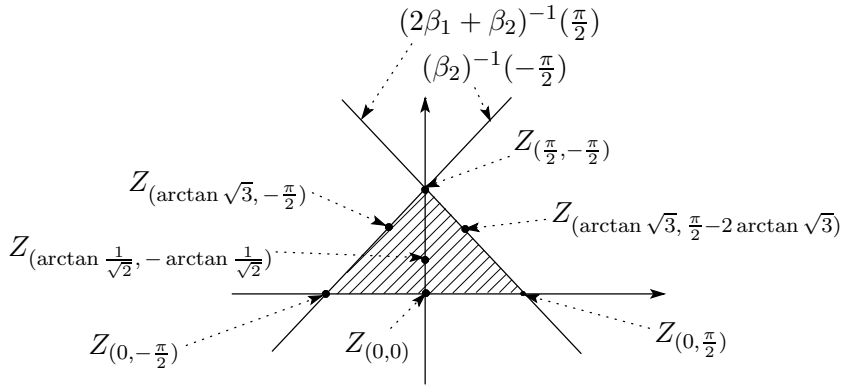
\centerline{\bf Figure 6.}

\newpage

\vspace{0.2truecm}

$$\begin{tabular}{|c|c|c|c|}
\hline
{\scriptsize$(a,b)$} & {\scriptsize $Z_{(a,b)}$} & 
{\scriptsize$M=U(1,1)^{\ast}({\rm Exp}\,Z_{(a,b)})$} & {\scriptsize ${\rm dim}\,M$}\\
\hline
{\scriptsize $(\frac{\pi}{2},0)$} & 
{\scriptsize as in Theorem F} & {\scriptsize one-point set} & $0$\\
\hline
{\scriptsize $(-\frac{\pi}{2},\pi)$} & 
{\scriptsize as in Theorem F} & {\scriptsize one-point set} & $0$\\
\hline
{\scriptsize $(0,0)$} & 
{\scriptsize as in Theorem F} & {\scriptsize totally geodesic} & $2$\\
\hline
{\scriptsize $(\frac{\pi}{6},0)$} & 
{\scriptsize as in Theorem C} & {\scriptsize not austere} & $3$\\
\hline
{\scriptsize $(-\frac{\pi}{6},\frac{\pi}{3})$} & 
{\scriptsize as in Theorem C} & {\scriptsize not austere} & $3$\\
\hline
{\scriptsize $(0,\frac{\pi}{2})$} & 
{\scriptsize as in Theorem F} & {\scriptsize totally geodesic} & $3$\\
\hline
{\scriptsize $(0,\arctan\sqrt2)$} & 
{\scriptsize not as in Theorems C$\sim$F} & {\scriptsize not austere} & $4$\\
\hline
\end{tabular}$$

\vspace{0.4truecm}

\centerline{$U(1,1)^{\ast}\curvearrowright Sp(2)/U(2)$}

\centerline{$({\rm dim}\,Sp(2)/U(2)=6)$}

\vspace{0.5truecm}

\centerline{{\bf Table 8.}}

\vspace{0.6truecm}

The positions of $Z_0$'s in Table 8 are as in Figure 7.  

\vspace{0.6truecm}

\centerline{
%\input{EXHEMF5-5.tex}
%WinTpicVersion3.08
\unitlength 0.1in
\begin{picture}( 47.3000, 19.8000)( -8.1000,-27.5000)
% VECTOR 2 0 3 0
% 2 2180 1940 3920 1940
% 
\special{pn 8}%
\special{pa 2180 1940}%
\special{pa 3920 1940}%
\special{fp}%
\special{sh 1}%
\special{pa 3920 1940}%
\special{pa 3854 1920}%
\special{pa 3868 1940}%
\special{pa 3854 1960}%
\special{pa 3920 1940}%
\special{fp}%
% VECTOR 2 0 3 0
% 4 2800 2680 2800 1060 2800 1060 2800 1060
% 
\special{pn 8}%
\special{pa 2800 2680}%
\special{pa 2800 1060}%
\special{fp}%
\special{sh 1}%
\special{pa 2800 1060}%
\special{pa 2780 1128}%
\special{pa 2800 1114}%
\special{pa 2820 1128}%
\special{pa 2800 1060}%
\special{fp}%
\special{pa 2800 1060}%
\special{pa 2800 1060}%
\special{fp}%
% STR 2 0 3 0
% 3 2240 1800 2240 1860 3 0
% $Z_{(0,0)}$
\put(22.4000,-18.6000){\makebox(0,0)[rb]{$Z_{(0,0)}$}}%
% STR 2 0 3 0
% 3 3800 1751 3800 1810 2 0
% $Z_{(0,\frac{\pi}{2})}$
\put(38.0000,-18.1000){\makebox(0,0)[lb]{$Z_{(0,\frac{\pi}{2})}$}}%
% STR 2 0 3 0
% 3 2660 896 2660 940 2 0
% $Z_{(0,\arctan\sqrt2)}$
\put(26.6000,-9.4000){\makebox(0,0)[lb]{$Z_{(0,\arctan\sqrt2)}$}}%
% STR 2 0 3 0
% 3 3810 2231 3810 2290 1 0
% $Z_{(-\frac{\pi}{2},\,\pi)}$
\put(38.1000,-22.9000){\makebox(0,0)[lt]{$Z_{(-\frac{\pi}{2},\,\pi)}$}}%
% STR 2 0 3 0
% 3 2620 1411 2620 1470 3 0
% $Z_{(\frac{\pi}{6},\,0)}$
\put(26.2000,-14.7000){\makebox(0,0)[rb]{$Z_{(\frac{\pi}{6},\,0)}$}}%
% STR 2 0 3 0
% 3 3580 1042 3580 1100 2 0
% $(\beta_1+\beta_2)^{-1}(\frac{\pi}{2})$
\put(35.8000,-11.0000){\makebox(0,0)[lb]{$(\beta_1+\beta_2)^{-1}(\frac{\pi}{2})$}}%
% LINE 2 0 3 0
% 2 2660 2080 3580 1221
% 
\special{pn 8}%
\special{pa 2660 2080}%
\special{pa 3580 1222}%
\special{fp}%
% LINE 2 0 3 0
% 2 2640 1780 3550 2650
% 
\special{pn 8}%
\special{pa 2640 1780}%
\special{pa 3550 2650}%
\special{fp}%
% LINE 3 0 3 0
% 40 2839 1959 2827 1934 2926 2043 2858 1905 3014 2127 2889 1877 3101 2211 2919 1849 3188 2295 2950 1821 3275 2380 2981 1792 3362 2464 3012 1764 3420 2490 3043 1736 3420 2400 3074 1707 3420 2310 3105 1679 3420 2220 3135 1651 3420 2130 3166 1623 3420 2040 3197 1594 3420 1950 3228 1566 3420 1860 3259 1538 3420 1770 3290 1509 3420 1680 3321 1481 3420 1590 3351 1453 3420 1500 3382 1425 3420 1410 3413 1396
% 
\special{pn 4}%
\special{pa 2840 1960}%
\special{pa 2828 1934}%
\special{fp}%
\special{pa 2926 2044}%
\special{pa 2858 1906}%
\special{fp}%
\special{pa 3014 2128}%
\special{pa 2890 1878}%
\special{fp}%
\special{pa 3102 2212}%
\special{pa 2920 1850}%
\special{fp}%
\special{pa 3188 2296}%
\special{pa 2950 1822}%
\special{fp}%
\special{pa 3276 2380}%
\special{pa 2982 1792}%
\special{fp}%
\special{pa 3362 2464}%
\special{pa 3012 1764}%
\special{fp}%
\special{pa 3420 2490}%
\special{pa 3044 1736}%
\special{fp}%
\special{pa 3420 2400}%
\special{pa 3074 1708}%
\special{fp}%
\special{pa 3420 2310}%
\special{pa 3106 1680}%
\special{fp}%
\special{pa 3420 2220}%
\special{pa 3136 1652}%
\special{fp}%
\special{pa 3420 2130}%
\special{pa 3166 1624}%
\special{fp}%
\special{pa 3420 2040}%
\special{pa 3198 1594}%
\special{fp}%
\special{pa 3420 1950}%
\special{pa 3228 1566}%
\special{fp}%
\special{pa 3420 1860}%
\special{pa 3260 1538}%
\special{fp}%
\special{pa 3420 1770}%
\special{pa 3290 1510}%
\special{fp}%
\special{pa 3420 1680}%
\special{pa 3322 1482}%
\special{fp}%
\special{pa 3420 1590}%
\special{pa 3352 1454}%
\special{fp}%
\special{pa 3420 1500}%
\special{pa 3382 1426}%
\special{fp}%
\special{pa 3420 1410}%
\special{pa 3414 1396}%
\special{fp}%
% DOT 0 0 3 0
% 2 2810 1940 2810 1940
% 
\special{pn 20}%
\special{sh 1}%
\special{ar 2810 1940 10 10 0  6.28318530717959E+0000}%
\special{sh 1}%
\special{ar 2810 1940 10 10 0  6.28318530717959E+0000}%
% DOT 0 0 3 0
% 2 3420 1370 3420 1370
% 
\special{pn 20}%
\special{sh 1}%
\special{ar 3420 1370 10 10 0  6.28318530717959E+0000}%
\special{sh 1}%
\special{ar 3420 1370 10 10 0  6.28318530717959E+0000}%
% DOT 0 0 3 0
% 2 3420 1940 3420 1940
% 
\special{pn 20}%
\special{sh 1}%
\special{ar 3420 1940 10 10 0  6.28318530717959E+0000}%
\special{sh 1}%
\special{ar 3420 1940 10 10 0  6.28318530717959E+0000}%
% DOT 0 0 3 0
% 2 3420 2530 3420 2530
% 
\special{pn 20}%
\special{sh 1}%
\special{ar 3420 2530 10 10 0  6.28318530717959E+0000}%
\special{sh 1}%
\special{ar 3420 2530 10 10 0  6.28318530717959E+0000}%
% DOT 0 0 3 0
% 2 3010 1750 3010 1750
% 
\special{pn 20}%
\special{sh 1}%
\special{ar 3010 1750 10 10 0  6.28318530717959E+0000}%
\special{sh 1}%
\special{ar 3010 1750 10 10 0  6.28318530717959E+0000}%
% DOT 0 0 3 0
% 2 3000 2130 3000 2130
% 
\special{pn 20}%
\special{sh 1}%
\special{ar 3000 2130 10 10 0  6.28318530717959E+0000}%
\special{sh 1}%
\special{ar 3000 2130 10 10 0  6.28318530717959E+0000}%
% DOT 0 0 3 0
% 2 3140 1940 3140 1940
% 
\special{pn 20}%
\special{sh 1}%
\special{ar 3140 1940 10 10 0  6.28318530717959E+0000}%
\special{sh 1}%
\special{ar 3140 1940 10 10 0  6.28318530717959E+0000}%
% LINE 2 0 3 0
% 2 3420 1200 3420 2690
% 
\special{pn 8}%
\special{pa 3420 1200}%
\special{pa 3420 2690}%
\special{fp}%
% VECTOR 2 2 3 0
% 2 3580 1100 3420 1270
% 
\special{pn 8}%
\special{pa 3580 1100}%
\special{pa 3420 1270}%
\special{dt 0.045}%
\special{sh 1}%
\special{pa 3420 1270}%
\special{pa 3480 1236}%
\special{pa 3458 1232}%
\special{pa 3452 1208}%
\special{pa 3420 1270}%
\special{fp}%
% VECTOR 2 2 3 0
% 2 3770 1780 3430 1930
% 
\special{pn 8}%
\special{pa 3770 1780}%
\special{pa 3430 1930}%
\special{dt 0.045}%
\special{sh 1}%
\special{pa 3430 1930}%
\special{pa 3500 1922}%
\special{pa 3480 1908}%
\special{pa 3484 1886}%
\special{pa 3430 1930}%
\special{fp}%
% VECTOR 2 2 3 0
% 2 3790 1480 3430 1380
% 
\special{pn 8}%
\special{pa 3790 1480}%
\special{pa 3430 1380}%
\special{dt 0.045}%
\special{sh 1}%
\special{pa 3430 1380}%
\special{pa 3490 1418}%
\special{pa 3482 1394}%
\special{pa 3500 1380}%
\special{pa 3430 1380}%
\special{fp}%
% VECTOR 2 2 3 0
% 2 2640 1450 3000 1750
% 
\special{pn 8}%
\special{pa 2640 1450}%
\special{pa 3000 1750}%
\special{dt 0.045}%
\special{sh 1}%
\special{pa 3000 1750}%
\special{pa 2962 1692}%
\special{pa 2960 1716}%
\special{pa 2936 1724}%
\special{pa 3000 1750}%
\special{fp}%
% VECTOR 2 2 3 0
% 2 2280 1810 2800 1940
% 
\special{pn 8}%
\special{pa 2280 1810}%
\special{pa 2800 1940}%
\special{dt 0.045}%
\special{sh 1}%
\special{pa 2800 1940}%
\special{pa 2740 1904}%
\special{pa 2748 1928}%
\special{pa 2730 1944}%
\special{pa 2800 1940}%
\special{fp}%
% VECTOR 2 2 3 0
% 2 2650 2350 2990 2140
% 
\special{pn 8}%
\special{pa 2650 2350}%
\special{pa 2990 2140}%
\special{dt 0.045}%
\special{sh 1}%
\special{pa 2990 2140}%
\special{pa 2924 2158}%
\special{pa 2946 2168}%
\special{pa 2944 2192}%
\special{pa 2990 2140}%
\special{fp}%
% VECTOR 2 2 3 0
% 2 3780 2390 3430 2530
% 
\special{pn 8}%
\special{pa 3780 2390}%
\special{pa 3430 2530}%
\special{dt 0.045}%
\special{sh 1}%
\special{pa 3430 2530}%
\special{pa 3500 2524}%
\special{pa 3480 2510}%
\special{pa 3484 2488}%
\special{pa 3430 2530}%
\special{fp}%
% VECTOR 2 2 3 0
% 2 2940 990 3140 1930
% 
\special{pn 8}%
\special{pa 2940 990}%
\special{pa 3140 1930}%
\special{dt 0.045}%
\special{sh 1}%
\special{pa 3140 1930}%
\special{pa 3146 1862}%
\special{pa 3130 1878}%
\special{pa 3108 1870}%
\special{pa 3140 1930}%
\special{fp}%
% STR 2 0 3 0
% 3 2610 2271 2610 2330 4 0
% $Z_{(-\frac{\pi}{6},\,\frac{\pi}{3})}$
\put(26.1000,-23.3000){\makebox(0,0)[rt]{$Z_{(-\frac{\pi}{6},\,\frac{\pi}{3})}$}}%
% VECTOR 2 2 3 0
% 2 2600 2750 3210 2340
% 
\special{pn 8}%
\special{pa 2600 2750}%
\special{pa 3210 2340}%
\special{dt 0.045}%
\special{sh 1}%
\special{pa 3210 2340}%
\special{pa 3144 2362}%
\special{pa 3166 2370}%
\special{pa 3166 2394}%
\special{pa 3210 2340}%
\special{fp}%
% VECTOR 2 2 3 0
% 2 2250 1030 3230 1540
% 
\special{pn 8}%
\special{pa 2250 1030}%
\special{pa 3230 1540}%
\special{dt 0.045}%
\special{sh 1}%
\special{pa 3230 1540}%
\special{pa 3180 1492}%
\special{pa 3184 1516}%
\special{pa 3162 1528}%
\special{pa 3230 1540}%
\special{fp}%
% STR 2 0 3 0
% 3 2550 2672 2550 2730 4 0
% $(2\beta_1+\beta_2)^{-1}(0)$
\put(25.5000,-27.3000){\makebox(0,0)[rt]{$(2\beta_1+\beta_2)^{-1}(0)$}}%
% STR 2 0 3 0
% 3 2200 992 2200 1050 3 0
% $\beta_2^{-1}(0)$
\put(22.0000,-10.5000){\makebox(0,0)[rb]{$\beta_2^{-1}(0)$}}%
% STR 2 0 3 0
% 3 3840 1511 3840 1570 2 0
% $Z_{(\frac{\pi}{2},0)}$
\put(38.4000,-15.7000){\makebox(0,0)[lb]{$Z_{(\frac{\pi}{2},0)}$}}%
\end{picture}%
\hspace{8.5truecm}}

\vspace{1truecm}

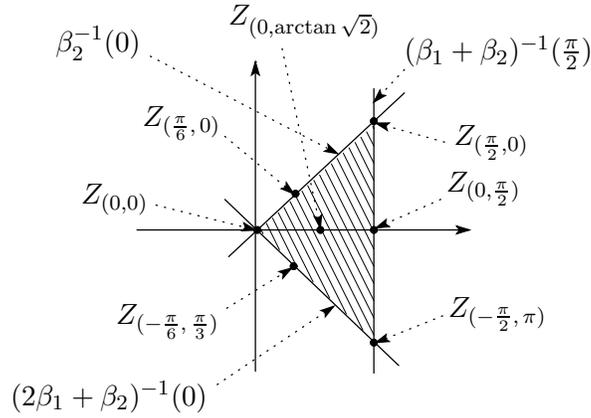
\centerline{\bf Figure 7.}

\newpage

\vspace{0.3truecm}

$$\begin{tabular}{|c|c|c|c|}
\hline
{\scriptsize$(a,b)$} & {\scriptsize $Z_{(a,b)}$} & 
{\scriptsize$M=Sp(2,{\Bbb R})^{\ast}({\rm Exp}\,Z_{(a,b)})$} & {\scriptsize ${\rm dim}\,M$}\\
\hline
{\scriptsize $(0,-\frac{\pi}{2})$} & 
{\scriptsize as in Theorem F} & {\scriptsize one-point set} & $0$\\
\hline
{\scriptsize $(0,\frac{\pi}{2})$} & 
{\scriptsize as in Theorem F} & {\scriptsize one-point set} & $0$\\
\hline
{\scriptsize $(\frac{\pi}{2},-\frac{\pi}{2})$} & 
{\scriptsize as in Theorem F} & {\scriptsize totally geodesic} & $4$\\
\hline
{\scriptsize $(0,0)$} & 
{\scriptsize as in Theorem F} & {\scriptsize totally geodesic} & $6$\\
\hline
{\scriptsize $(\arctan\sqrt 3,-\frac{\pi}{2})$} & 
{\scriptsize not as in Theorems C$\sim$F} & {\scriptsize not austere} & $6$\\
\hline
{\scriptsize $(\arctan\sqrt 3,\frac{\pi}{2}-2\arctan\sqrt 3)$} & 
{\scriptsize not as in Theorems C$\sim$F} & {\scriptsize not austere} & $6$\\
\hline
{\scriptsize $(\arctan\frac{1}{\sqrt 2},-\arctan\frac{1}{\sqrt2})$} & 
{\scriptsize not as in Theorems C$\sim$F} & {\scriptsize not austere} & $8$\\
\hline
\end{tabular}$$

\vspace{0.5truecm}

\centerline{$Sp(2,{\Bbb R})^{\ast}\curvearrowright(Sp(2)\times Sp(2))/Sp(2)$}

\centerline{$({\rm dim}\,(Sp(2)\times Sp(2))/Sp(2)=10)$}

\vspace{0.5truecm}

\centerline{{\bf Table 9.}}

%\newpage

\vspace{0.5truecm}

The positions of $Z_0$'s in Table 9 are as in Figure 6.  

\vspace{0.5truecm}

$$\begin{tabular}{|c|c|c|c|}
\hline
{\scriptsize$(a,b)$} & {\scriptsize $Z_{(a,b)}$} & 
{\scriptsize$M=Sp(1,1)^{\ast}({\rm Exp}\,Z_{(a,b)})$} & {\scriptsize ${\rm dim}\,M$}\\
\hline
{\scriptsize $(\frac{\pi}{2},0)$} & 
{\scriptsize as in Theorem F} & {\scriptsize one-point set} & $0$\\
\hline
{\scriptsize $(-\frac{\pi}{2},\pi)$} & 
{\scriptsize as in Theorem F} & {\scriptsize one-point set} & $0$\\
\hline
{\scriptsize $(0,0)$} & 
{\scriptsize as in Theorem F} & {\scriptsize totally geodesic} & $4$\\
\hline
{\scriptsize $(\frac{\pi}{6},0)$} & 
{\scriptsize as in Theorem C} & {\scriptsize not austere} & $6$\\
\hline
{\scriptsize $(-\frac{\pi}{6},\frac{\pi}{3})$} & 
{\scriptsize as in Theorem C} & {\scriptsize not austere} & $6$\\
\hline
{\scriptsize $(0,\frac{\pi}{2})$} & 
{\scriptsize as in Theorem F} & {\scriptsize totally geodesic} & $6$\\
\hline
{\scriptsize $(0,\arctan\sqrt 2)$} & 
{\scriptsize not as in Theorems C$\sim$F} & {\scriptsize not austere} & $8$\\
\hline
\end{tabular}$$

\vspace{0.5truecm}

\centerline{$Sp(1,1)^{\ast}\curvearrowright(Sp(2)\times Sp(2))/Sp(2)$}

\centerline{$({\rm dim}\,(Sp(2)\times Sp(2))/Sp(2)=10)$}

\vspace{0.5truecm}

\centerline{{\bf Table 10.}}

\vspace{0.5truecm}

The positions of $Z_0$'s in Table 10 are as in Figure 7.  

\newpage

%\vspace{0.5truecm}

$$\begin{tabular}{|c|c|c|c|}
\hline
{\scriptsize$(a,b)$} & {\scriptsize $Z_{(a,b)}$} & 
{\scriptsize$M=(SO^{\ast}(10)\cdot U(1))^{\ast}({\rm Exp}\,Z_{(a,b)})$} & 
{\scriptsize ${\rm dim}\,M$}\\
\hline
{\scriptsize $(0,0)$} & 
{\scriptsize as in Theorem F} & {\scriptsize totally geodesic} & $20$\\
\hline
{\scriptsize $(0,\frac{\pi}{2})$} & 
{\scriptsize as in Theorem F} & {\scriptsize one-point set} & $0$\\
\hline
{\scriptsize $(\frac{\pi}{2},-\frac{\pi}{2})$} & 
{\scriptsize as in Theorem F} & {\scriptsize totally geodesic} & $17$\\
\hline
{\scriptsize $(0,a_1)$} & 
{\scriptsize not as in Theorems C$\sim$F} & {\scriptsize not austere} & $21$\\
\hline
{\scriptsize $(a_2,-a_2)$} & 
{\scriptsize not as in Theorems C$\sim$F} & {\scriptsize not austere} & $29$\\
\hline
{\scriptsize $(a_3,\frac{\pi}{2}-2a_3)$} & 
{\scriptsize not as in Theorems C$\sim$F} & {\scriptsize not austere} & $25$\\
\hline
{\scriptsize $(a_4,b)$} & 
{\scriptsize not as in Theorems C$\sim$F} & {\scriptsize not austere} & $30$\\
\hline
\end{tabular}$$

\vspace{0.5truecm}

\centerline{$(SO^{\ast}(10)\cdot U(1))^{\ast}\curvearrowright 
E_6/Spin(10)\cdot U(1)$}

\centerline{$({\rm dim}\,E_6/Spin(10)\cdot U(1)=32)$}

\vspace{0.5truecm}

\centerline{{\bf Table 11.}}

\vspace{0.5truecm}

The positions of $Z_0$'s in Table 11 are as in Figure 8.  
The numbers $a_i$ ($i=1,2,3,4$) and $b$ in Table 11 are real numbers such that 
$a_i,b\not\equiv\frac{\pi}{6},\frac{\pi}{3},\frac{\pi}{4},\frac{3\pi}{4}\,\,
({\rm mod}\,\pi)$.  

\vspace{0.5truecm}

\centerline{
%\input{EXHEMF6.5.tex}
%WinTpicVersion3.08
\unitlength 0.1in
\begin{picture}( 49.0000, 23.2000)( -2.9000,-27.0000)
% VECTOR 2 0 3 0
% 2 2350 2110 4610 2110
% 
\special{pn 8}%
\special{pa 2350 2110}%
\special{pa 4610 2110}%
\special{fp}%
\special{sh 1}%
\special{pa 4610 2110}%
\special{pa 4544 2090}%
\special{pa 4558 2110}%
\special{pa 4544 2130}%
\special{pa 4610 2110}%
\special{fp}%
% VECTOR 2 0 3 0
% 4 3200 2600 3200 559 3200 559 3200 559
% 
\special{pn 8}%
\special{pa 3200 2600}%
\special{pa 3200 560}%
\special{fp}%
\special{sh 1}%
\special{pa 3200 560}%
\special{pa 3180 626}%
\special{pa 3200 612}%
\special{pa 3220 626}%
\special{pa 3200 560}%
\special{fp}%
\special{pa 3200 560}%
\special{pa 3200 560}%
\special{fp}%
% LINE 2 0 3 0
% 2 4210 1570 2386 1568
% 
\special{pn 8}%
\special{pa 4210 1570}%
\special{pa 2386 1568}%
\special{fp}%
% DOT 0 0 3 0
% 2 3250 1770 3250 1770
% 
\special{pn 20}%
\special{sh 1}%
\special{ar 3250 1770 10 10 0  6.28318530717959E+0000}%
\special{sh 1}%
\special{ar 3250 1770 10 10 0  6.28318530717959E+0000}%
% STR 2 0 3 0
% 3 4210 1292 4210 1350 2 0
% $(2\beta_1+\beta_2)^{-1}(\frac{\pi}{2})$
\put(42.1000,-13.5000){\makebox(0,0)[lb]{$(2\beta_1+\beta_2)^{-1}(\frac{\pi}{2})$}}%
% STR 2 0 3 0
% 3 3280 2592 3280 2660 1 0
% $(\beta_1)^{-1}(0)$
\put(32.8000,-26.6000){\makebox(0,0)[lt]{$(\beta_1)^{-1}(0)$}}%
% STR 2 0 3 0
% 3 2950 872 2950 940 3 0
% $Z_{(\frac{\pi}{2},-\frac{\pi}{2})}$
\put(29.5000,-9.4000){\makebox(0,0)[rb]{$Z_{(\frac{\pi}{2},-\frac{\pi}{2})}$}}%
% DOT 0 0 3 0
% 2 3250 1570 3250 1570
% 
\special{pn 20}%
\special{sh 1}%
\special{ar 3250 1570 10 10 0  6.28318530717959E+0000}%
\special{sh 1}%
\special{ar 3250 1570 10 10 0  6.28318530717959E+0000}%
% DOT 2 0 3 0
% 2 3721 1568 3721 1568
% 
\special{pn 8}%
\special{sh 1}%
\special{ar 3722 1568 10 10 0  6.28318530717959E+0000}%
\special{sh 1}%
\special{ar 3722 1568 10 10 0  6.28318530717959E+0000}%
% STR 2 0 3 0
% 3 4200 1621 4200 1680 1 0
% $Z_{(0,\frac{\pi}{2})}$
\put(42.0000,-16.8000){\makebox(0,0)[lt]{$Z_{(0,\frac{\pi}{2})}$}}%
% STR 2 0 3 0
% 3 2760 312 2760 380 4 0
% $(\beta_2)^{-1}(0)$
\put(27.6000,-3.8000){\makebox(0,0)[rt]{$(\beta_2)^{-1}(0)$}}%
% STR 2 0 3 0
% 3 2480 2220 2480 2320 4 0
% $Z_{(a_4,b)}$
\put(24.8000,-23.2000){\makebox(0,0)[rt]{$Z_{(a_4,b)}$}}%
% LINE 2 0 3 0
% 2 4020 1250 2770 2578
% 
\special{pn 8}%
\special{pa 4020 1250}%
\special{pa 2770 2578}%
\special{fp}%
% DOT 0 0 3 0
% 2 3210 2110 3210 2110
% 
\special{pn 20}%
\special{sh 1}%
\special{ar 3210 2110 10 10 0  6.28318530717959E+0000}%
\special{sh 1}%
\special{ar 3210 2110 10 10 0  6.28318530717959E+0000}%
% DOT 0 0 3 0
% 2 3520 1780 3520 1780
% 
\special{pn 20}%
\special{sh 1}%
\special{ar 3520 1780 10 10 0  6.28318530717959E+0000}%
\special{sh 1}%
\special{ar 3520 1780 10 10 0  6.28318530717959E+0000}%
% STR 2 0 3 0
% 3 3850 990 3850 1090 2 0
% $Z_{(a_3,\frac{\pi}{2}-2a_3)}$
\put(38.5000,-10.9000){\makebox(0,0)[lb]{$Z_{(a_3,\frac{\pi}{2}-2a_3)}$}}%
% STR 2 0 3 0
% 3 4210 1800 4210 1900 1 0
% $Z_{(0,a_1)}$
\put(42.1000,-19.0000){\makebox(0,0)[lt]{$Z_{(0,a_1)}$}}%
% VECTOR 2 2 3 0
% 2 3250 2700 2860 2490
% 
\special{pn 8}%
\special{pa 3250 2700}%
\special{pa 2860 2490}%
\special{dt 0.045}%
\special{sh 1}%
\special{pa 2860 2490}%
\special{pa 2910 2540}%
\special{pa 2908 2516}%
\special{pa 2928 2504}%
\special{pa 2860 2490}%
\special{fp}%
% STR 2 0 3 0
% 3 3900 2250 3900 2300 1 0
% $Z_{(0,0)}$
\put(39.0000,-23.0000){\makebox(0,0)[lt]{$Z_{(0,0)}$}}%
% VECTOR 2 2 3 0
% 2 4170 1980 3530 1790
% 
\special{pn 8}%
\special{pa 4170 1980}%
\special{pa 3530 1790}%
\special{dt 0.045}%
\special{sh 1}%
\special{pa 3530 1790}%
\special{pa 3588 1828}%
\special{pa 3582 1806}%
\special{pa 3600 1790}%
\special{pa 3530 1790}%
\special{fp}%
% DOT 0 0 3 0
% 2 3730 1570 3730 1570
% 
\special{pn 20}%
\special{sh 1}%
\special{ar 3730 1570 10 10 0  6.28318530717959E+0000}%
\special{sh 1}%
\special{ar 3730 1570 10 10 0  6.28318530717959E+0000}%
% VECTOR 2 2 3 0
% 2 2420 1890 2900 1800
% 
\special{pn 8}%
\special{pa 2420 1890}%
\special{pa 2900 1800}%
\special{dt 0.045}%
\special{sh 1}%
\special{pa 2900 1800}%
\special{pa 2832 1794}%
\special{pa 2848 1810}%
\special{pa 2838 1832}%
\special{pa 2900 1800}%
\special{fp}%
% VECTOR 2 2 3 0
% 2 4180 1760 3740 1580
% 
\special{pn 8}%
\special{pa 4180 1760}%
\special{pa 3740 1580}%
\special{dt 0.045}%
\special{sh 1}%
\special{pa 3740 1580}%
\special{pa 3794 1624}%
\special{pa 3790 1600}%
\special{pa 3810 1588}%
\special{pa 3740 1580}%
\special{fp}%
% VECTOR 2 2 3 0
% 2 4310 1370 4100 1570
% 
\special{pn 8}%
\special{pa 4310 1370}%
\special{pa 4100 1570}%
\special{dt 0.045}%
\special{sh 1}%
\special{pa 4100 1570}%
\special{pa 4162 1540}%
\special{pa 4140 1534}%
\special{pa 4134 1510}%
\special{pa 4100 1570}%
\special{fp}%
% LINE 2 0 3 0
% 2 2390 1240 3640 2568
% 
\special{pn 8}%
\special{pa 2390 1240}%
\special{pa 3640 2568}%
\special{fp}%
% DOT 0 0 3 0
% 2 2690 1570 2690 1570
% 
\special{pn 20}%
\special{sh 1}%
\special{ar 2690 1570 10 10 0  6.28318530717959E+0000}%
\special{sh 1}%
\special{ar 2690 1570 10 10 0  6.28318530717959E+0000}%
% DOT 0 0 3 0
% 2 2920 1800 2920 1800
% 
\special{pn 20}%
\special{sh 1}%
\special{ar 2920 1800 10 10 0  6.28318530717959E+0000}%
\special{sh 1}%
\special{ar 2920 1800 10 10 0  6.28318530717959E+0000}%
% VECTOR 2 2 3 0
% 2 2590 980 2700 1570
% 
\special{pn 8}%
\special{pa 2590 980}%
\special{pa 2700 1570}%
\special{dt 0.045}%
\special{sh 1}%
\special{pa 2700 1570}%
\special{pa 2708 1502}%
\special{pa 2690 1518}%
\special{pa 2668 1508}%
\special{pa 2700 1570}%
\special{fp}%
% VECTOR 2 2 3 0
% 2 2790 510 3200 730
% 
\special{pn 8}%
\special{pa 2790 510}%
\special{pa 3200 730}%
\special{dt 0.045}%
\special{sh 1}%
\special{pa 3200 730}%
\special{pa 3152 682}%
\special{pa 3154 706}%
\special{pa 3132 716}%
\special{pa 3200 730}%
\special{fp}%
% STR 2 0 3 0
% 3 2370 1740 2370 1840 4 0
% $Z_{(a_2,-a_2)}$
\put(23.7000,-18.4000){\makebox(0,0)[rt]{$Z_{(a_2,-a_2)}$}}%
% VECTOR 2 2 3 0
% 2 2490 2350 3250 1780
% 
\special{pn 8}%
\special{pa 2490 2350}%
\special{pa 3250 1780}%
\special{dt 0.045}%
\special{sh 1}%
\special{pa 3250 1780}%
\special{pa 3186 1804}%
\special{pa 3208 1812}%
\special{pa 3210 1836}%
\special{pa 3250 1780}%
\special{fp}%
% VECTOR 2 2 3 0
% 2 3870 2370 3210 2120
% 
\special{pn 8}%
\special{pa 3870 2370}%
\special{pa 3210 2120}%
\special{dt 0.045}%
\special{sh 1}%
\special{pa 3210 2120}%
\special{pa 3266 2162}%
\special{pa 3260 2140}%
\special{pa 3280 2126}%
\special{pa 3210 2120}%
\special{fp}%
% VECTOR 2 2 3 0
% 2 3830 1070 3260 1560
% 
\special{pn 8}%
\special{pa 3830 1070}%
\special{pa 3260 1560}%
\special{dt 0.045}%
\special{sh 1}%
\special{pa 3260 1560}%
\special{pa 3324 1532}%
\special{pa 3300 1526}%
\special{pa 3298 1502}%
\special{pa 3260 1560}%
\special{fp}%
% LINE 3 0 3 0
% 34 3225 2085 2720 1580 3254 2054 2780 1580 3284 2024 2840 1580 3313 1993 2900 1580 3342 1962 2960 1580 3372 1932 3020 1580 3401 1901 3080 1580 3431 1871 3140 1580 3460 1840 3200 1580 3489 1809 3260 1580 3519 1779 3320 1580 3548 1748 3380 1580 3578 1718 3440 1580 3607 1687 3500 1580 3636 1656 3560 1580 3666 1626 3620 1580 3695 1595 3680 1580
% 
\special{pn 4}%
\special{pa 3226 2086}%
\special{pa 2720 1580}%
\special{fp}%
\special{pa 3254 2054}%
\special{pa 2780 1580}%
\special{fp}%
\special{pa 3284 2024}%
\special{pa 2840 1580}%
\special{fp}%
\special{pa 3314 1994}%
\special{pa 2900 1580}%
\special{fp}%
\special{pa 3342 1962}%
\special{pa 2960 1580}%
\special{fp}%
\special{pa 3372 1932}%
\special{pa 3020 1580}%
\special{fp}%
\special{pa 3402 1902}%
\special{pa 3080 1580}%
\special{fp}%
\special{pa 3432 1872}%
\special{pa 3140 1580}%
\special{fp}%
\special{pa 3460 1840}%
\special{pa 3200 1580}%
\special{fp}%
\special{pa 3490 1810}%
\special{pa 3260 1580}%
\special{fp}%
\special{pa 3520 1780}%
\special{pa 3320 1580}%
\special{fp}%
\special{pa 3548 1748}%
\special{pa 3380 1580}%
\special{fp}%
\special{pa 3578 1718}%
\special{pa 3440 1580}%
\special{fp}%
\special{pa 3608 1688}%
\special{pa 3500 1580}%
\special{fp}%
\special{pa 3636 1656}%
\special{pa 3560 1580}%
\special{fp}%
\special{pa 3666 1626}%
\special{pa 3620 1580}%
\special{fp}%
\special{pa 3696 1596}%
\special{pa 3680 1580}%
\special{fp}%
% STR 2 0 3 0
% 3 2230 1272 2230 1340 4 0
% $(\beta_1+\beta_2)^{-1}(0)$
\put(22.3000,-13.4000){\makebox(0,0)[rt]{$(\beta_1+\beta_2)^{-1}(0)$}}%
% VECTOR 2 2 3 0
% 2 2250 1400 2450 1310
% 
\special{pn 8}%
\special{pa 2250 1400}%
\special{pa 2450 1310}%
\special{dt 0.045}%
\special{sh 1}%
\special{pa 2450 1310}%
\special{pa 2382 1320}%
\special{pa 2402 1332}%
\special{pa 2398 1356}%
\special{pa 2450 1310}%
\special{fp}%
\end{picture}%
\hspace{6.2truecm}}

\vspace{0.5truecm}

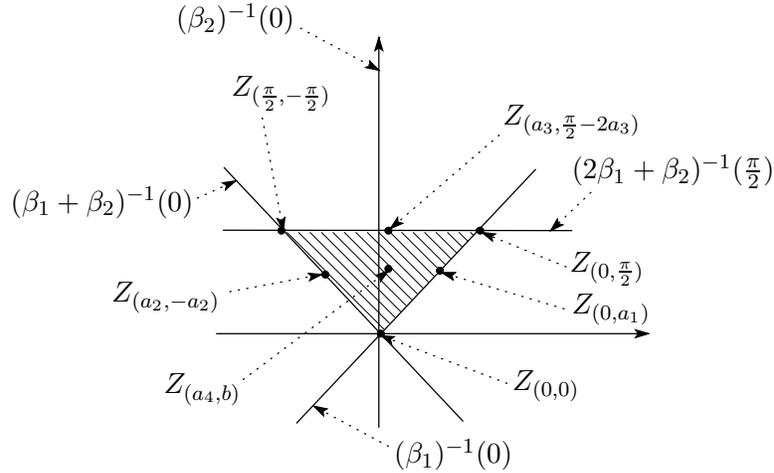
\centerline{\bf Figure 8.}

\newpage

%\vspace{0.3truecm}

$$\begin{tabular}{|c|c|c|c|}
\hline
{\scriptsize$(a,b)$} & {\scriptsize $Z_{(a,b)}$} & 
{\scriptsize$M=(F_4^{-20})^{\ast}({\rm Exp}\,Z_{(a,b)})$} & {\scriptsize ${\rm dim}\,M$}\\
\hline
{\scriptsize $(0,-\frac{\pi}{2})$} & 
{\scriptsize as in Theorem F} & {\scriptsize one-point set} & $0$\\
\hline
{\scriptsize $(0,\frac{\pi}{2})$} & 
{\scriptsize as in Theorem F} & {\scriptsize one-point set} & $0$\\
\hline
{\scriptsize $(\pi,-\frac{\pi}{2})$} & 
{\scriptsize as in Theorem F} & {\scriptsize one-point set} & $0$\\
\hline
{\scriptsize $(0,0)$} & 
{\scriptsize as in Theorem F} & {\scriptsize totally geodesic} & $16$\\
\hline
{\scriptsize $(\frac{\pi}{2},0)$} & 
{\scriptsize as in Theorem F} & {\scriptsize totally geodesic} & $16$\\
\hline
{\scriptsize $(\frac{\pi}{2},-\frac{\pi}{2})$} & 
{\scriptsize as in Theorem F} & {\scriptsize totally geodesic} & $16$\\
\hline
{\scriptsize $(\frac{\pi}{3},-\frac{\pi}{6})$} & 
{\scriptsize as in Theorem C} & {\scriptsize not austere} & $24$\\
\hline
\end{tabular}$$

\vspace{0.25truecm}

\centerline{$(F_4^{-20})^{\ast}\curvearrowright E_6/F_4$}

\centerline{$({\rm dim}\,E_6/F_4=26)$}

\vspace{0.25truecm}

\centerline{{\bf Table 12.}}

\vspace{0.35truecm}

The positions of $Z_0$'s in Table 12 are as in Figure 4.  

\vspace{0.25truecm}

%\newpage

$$\begin{tabular}{|c|c|c|c|}
\hline
{\scriptsize$(a,b)$} & {\scriptsize $Z_{(a,b)}$} & 
{\scriptsize$M=(SL(2,{\Bbb R})\times SL(2,{\Bbb R}))^{\ast}({\rm Exp}\,Z_{(a,b)})$} & 
{\scriptsize ${\rm dim}\,M$}\\
\hline
{\scriptsize $(0,-\frac{\pi}{2})$} & 
{\scriptsize as in Theorem F} & {\scriptsize one-point set} & $0$\\
\hline
{\scriptsize $(0,\frac{\pi}{2})$} & 
{\scriptsize as in Theorem F} & {\scriptsize one-point set} & $0$\\
\hline
{\scriptsize $(\frac{\pi}{2},-\frac{\pi}{2})$} & 
{\scriptsize as in Theorem F} & {\scriptsize totally geodesic} & $4$\\
\hline
{\scriptsize $(\frac{\pi}{3},-\frac{\pi}{2})$} & 
{\scriptsize as in Theorem C} & {\scriptsize not austere} & $3$\\
\hline
{\scriptsize $(\arctan\sqrt 5,\frac{\pi}{2}-2\arctan\sqrt 5)$} & 
{\scriptsize not as in Theorems C$\sim$F} & {\scriptsize not austere} & $5$\\
\hline
{\scriptsize $(a_4,b_2)$} & 
{\scriptsize not as in Theorems C$\sim$F} & {\scriptsize not austere} & $6$\\
\hline
\end{tabular}$$

\vspace{0.2truecm}

\centerline{$(SL(2,{\Bbb R})\times SL(2,{\Bbb R}))^{\ast}\curvearrowright G_2/SO(4)$}

\centerline{$({\rm dim}\,G_2/SO(4)=8)$}

\vspace{0.25truecm}

\centerline{{\bf Table 13.}}

\vspace{0.25truecm}

The positions of $Z_0$'s in Table 13 are as in Figure 9.  
The numbers $a_4$ and $b_2$ in Table 13 are real numbers such that 
$a_4,b_2\not\equiv\frac{\pi}{6},\frac{\pi}{3},\frac{\pi}{4},
\frac{3\pi}{4}$\newline
$({\rm mod}\,\pi)$.  

\vspace{0.2truecm}

\centerline{
%\input{EXHEMF7.tex}
%WinTpicVersion3.08
\unitlength 0.1in
\begin{picture}( 46.2200, 21.8000)( -9.8000,-24.8000)
% VECTOR 2 0 3 0
% 2 1372 1896 3642 1896
% 
\special{pn 8}%
\special{pa 1372 1896}%
\special{pa 3642 1896}%
\special{fp}%
\special{sh 1}%
\special{pa 3642 1896}%
\special{pa 3576 1876}%
\special{pa 3590 1896}%
\special{pa 3576 1916}%
\special{pa 3642 1896}%
\special{fp}%
% VECTOR 2 0 3 0
% 4 2590 2480 2590 643 2590 643 2590 643
% 
\special{pn 8}%
\special{pa 2590 2480}%
\special{pa 2590 644}%
\special{fp}%
\special{sh 1}%
\special{pa 2590 644}%
\special{pa 2570 710}%
\special{pa 2590 696}%
\special{pa 2610 710}%
\special{pa 2590 644}%
\special{fp}%
\special{pa 2590 644}%
\special{pa 2590 644}%
\special{fp}%
% LINE 2 0 3 0
% 2 1910 2070 3530 1145
% 
\special{pn 8}%
\special{pa 1910 2070}%
\special{pa 3530 1146}%
\special{fp}%
% LINE 2 0 3 0
% 2 3181 2285 2274 655
% 
\special{pn 8}%
\special{pa 3182 2286}%
\special{pa 2274 656}%
\special{fp}%
% DOT 0 0 3 0
% 2 2592 1896 2592 1896
% 
\special{pn 20}%
\special{sh 1}%
\special{ar 2592 1896 10 10 0  6.28318530717959E+0000}%
\special{sh 1}%
\special{ar 2592 1896 10 10 0  6.28318530717959E+0000}%
% DOT 0 0 3 0
% 2 2789 1566 2789 1566
% 
\special{pn 20}%
\special{sh 1}%
\special{ar 2790 1566 10 10 0  6.28318530717959E+0000}%
\special{sh 1}%
\special{ar 2790 1566 10 10 0  6.28318530717959E+0000}%
% VECTOR 2 2 3 0
% 2 3424 1471 2797 1557
% 
\special{pn 8}%
\special{pa 3424 1472}%
\special{pa 2798 1558}%
\special{dt 0.045}%
\special{sh 1}%
\special{pa 2798 1558}%
\special{pa 2866 1568}%
\special{pa 2850 1550}%
\special{pa 2860 1528}%
\special{pa 2798 1558}%
\special{fp}%
% VECTOR 2 2 3 0
% 2 2830 410 2296 726
% 
\special{pn 8}%
\special{pa 2830 410}%
\special{pa 2296 726}%
\special{dt 0.045}%
\special{sh 1}%
\special{pa 2296 726}%
\special{pa 2364 710}%
\special{pa 2342 700}%
\special{pa 2344 676}%
\special{pa 2296 726}%
\special{fp}%
% STR 2 0 3 0
% 3 2860 385 2860 470 2 0
% $(2\beta_1+\beta_2)^{-1}(\frac{\pi}{2})$
\put(28.6000,-4.7000){\makebox(0,0)[lb]{$(2\beta_1+\beta_2)^{-1}(\frac{\pi}{2})$}}%
% STR 2 0 3 0
% 3 3330 855 3330 940 2 0
% $(\beta_2)^{-1}(-\frac{\pi}{2})$
\put(33.3000,-9.4000){\makebox(0,0)[lb]{$(\beta_2)^{-1}(-\frac{\pi}{2})$}}%
% STR 2 0 3 0
% 3 1720 2225 1720 2310 4 0
% $(\beta_1)^{-1}(0)$
\put(17.2000,-23.1000){\makebox(0,0)[rt]{$(\beta_1)^{-1}(0)$}}%
% STR 2 0 3 0
% 3 2920 2385 2920 2470 1 0
% $Z_{(a_4,b_2)}$
\put(29.2000,-24.7000){\makebox(0,0)[lt]{$Z_{(a_4,b_2)}$}}%
% STR 2 0 3 0
% 3 2490 2055 2490 2140 4 0
% $Z_{(0,0)}$
\put(24.9000,-21.4000){\makebox(0,0)[rt]{$Z_{(0,0)}$}}%
% STR 2 0 3 0
% 3 2260 1385 2260 1470 3 0
% $Z_{(\frac{\pi}{3},-\frac{\pi}{2})}$
\put(22.6000,-14.7000){\makebox(0,0)[rb]{$Z_{(\frac{\pi}{3},-\frac{\pi}{2})}$}}%
% DOT 0 0 3 0
% 2 2210 1900 2210 1900
% 
\special{pn 20}%
\special{sh 1}%
\special{ar 2210 1900 10 10 0  6.28318530717959E+0000}%
\special{sh 1}%
\special{ar 2210 1900 10 10 0  6.28318530717959E+0000}%
% DOT 2 0 3 0
% 2 2960 1900 2960 1900
% 
\special{pn 8}%
\special{sh 1}%
\special{ar 2960 1900 10 10 0  6.28318530717959E+0000}%
\special{sh 1}%
\special{ar 2960 1900 10 10 0  6.28318530717959E+0000}%
% VECTOR 2 2 3 0
% 2 3320 1710 2860 1680
% 
\special{pn 8}%
\special{pa 3320 1710}%
\special{pa 2860 1680}%
\special{dt 0.045}%
\special{sh 1}%
\special{pa 2860 1680}%
\special{pa 2926 1704}%
\special{pa 2914 1684}%
\special{pa 2928 1664}%
\special{pa 2860 1680}%
\special{fp}%
% VECTOR 2 2 3 0
% 2 1860 1740 2204 1894
% 
\special{pn 8}%
\special{pa 1860 1740}%
\special{pa 2204 1894}%
\special{dt 0.045}%
\special{sh 1}%
\special{pa 2204 1894}%
\special{pa 2152 1850}%
\special{pa 2156 1872}%
\special{pa 2136 1886}%
\special{pa 2204 1894}%
\special{fp}%
% STR 2 0 3 0
% 3 1840 1705 1840 1790 3 0
% $Z_{(0,-\frac{\pi}{2})}$
\put(18.4000,-17.9000){\makebox(0,0)[rb]{$Z_{(0,-\frac{\pi}{2})}$}}%
% VECTOR 2 2 3 0
% 2 3440 2120 2970 1910
% 
\special{pn 8}%
\special{pa 3440 2120}%
\special{pa 2970 1910}%
\special{dt 0.045}%
\special{sh 1}%
\special{pa 2970 1910}%
\special{pa 3024 1956}%
\special{pa 3020 1932}%
\special{pa 3040 1920}%
\special{pa 2970 1910}%
\special{fp}%
% STR 2 0 3 0
% 3 3480 2015 3480 2100 1 0
% $Z_{(0,\frac{\pi}{2})}$
\put(34.8000,-21.0000){\makebox(0,0)[lt]{$Z_{(0,\frac{\pi}{2})}$}}%
% STR 2 0 3 0
% 3 3440 1305 3440 1390 1 0
% $Z_{(\frac{\pi}{2},-\frac{\pi}{2})}$
\put(34.4000,-13.9000){\makebox(0,0)[lt]{$Z_{(\frac{\pi}{2},-\frac{\pi}{2})}$}}%
% DOT 0 0 3 0
% 2 2590 1680 2590 1680
% 
\special{pn 20}%
\special{sh 1}%
\special{ar 2590 1680 10 10 0  6.28318530717959E+0000}%
\special{sh 1}%
\special{ar 2590 1680 10 10 0  6.28318530717959E+0000}%
% STR 2 0 3 0
% 3 3360 1565 3360 1650 1 0
% $Z_{(\arctan\sqrt 5,\frac{\pi}{2}-2\arctan\sqrt 5)}$
\put(33.6000,-16.5000){\makebox(0,0)[lt]{$Z_{(\arctan\sqrt 5,\frac{\pi}{2}-2\arctan\sqrt 5)}$}}%
% VECTOR 2 2 3 0
% 2 2350 2110 2580 1910
% 
\special{pn 8}%
\special{pa 2350 2110}%
\special{pa 2580 1910}%
\special{dt 0.045}%
\special{sh 1}%
\special{pa 2580 1910}%
\special{pa 2518 1940}%
\special{pa 2540 1946}%
\special{pa 2544 1970}%
\special{pa 2580 1910}%
\special{fp}%
% VECTOR 2 2 3 0
% 2 1470 2260 1580 1900
% 
\special{pn 8}%
\special{pa 1470 2260}%
\special{pa 1580 1900}%
\special{dt 0.045}%
\special{sh 1}%
\special{pa 1580 1900}%
\special{pa 1542 1958}%
\special{pa 1564 1952}%
\special{pa 1580 1970}%
\special{pa 1580 1900}%
\special{fp}%
% DOT 0 0 3 0
% 2 2850 1690 2850 1690
% 
\special{pn 20}%
\special{sh 1}%
\special{ar 2850 1690 10 10 0  6.28318530717959E+0000}%
\special{sh 1}%
\special{ar 2850 1690 10 10 0  6.28318530717959E+0000}%
% VECTOR 2 2 3 0
% 2 2270 1450 2570 1680
% 
\special{pn 8}%
\special{pa 2270 1450}%
\special{pa 2570 1680}%
\special{dt 0.045}%
\special{sh 1}%
\special{pa 2570 1680}%
\special{pa 2530 1624}%
\special{pa 2528 1648}%
\special{pa 2506 1656}%
\special{pa 2570 1680}%
\special{fp}%
% DOT 0 0 3 0
% 2 2690 1770 2690 1770
% 
\special{pn 20}%
\special{sh 1}%
\special{ar 2690 1770 10 10 0  6.28318530717959E+0000}%
\special{sh 1}%
\special{ar 2690 1770 10 10 0  6.28318530717959E+0000}%
% VECTOR 2 2 3 0
% 2 3440 970 3340 1250
% 
\special{pn 8}%
\special{pa 3440 970}%
\special{pa 3340 1250}%
\special{dt 0.045}%
\special{sh 1}%
\special{pa 3340 1250}%
\special{pa 3382 1194}%
\special{pa 3358 1200}%
\special{pa 3344 1180}%
\special{pa 3340 1250}%
\special{fp}%
% VECTOR 2 2 3 0
% 2 3010 2430 2690 1770
% 
\special{pn 8}%
\special{pa 3010 2430}%
\special{pa 2690 1770}%
\special{dt 0.045}%
\special{sh 1}%
\special{pa 2690 1770}%
\special{pa 2702 1840}%
\special{pa 2714 1818}%
\special{pa 2738 1822}%
\special{pa 2690 1770}%
\special{fp}%
% LINE 3 0 3 0
% 24 2240 1900 2277 1863 2300 1900 2411 1789 2360 1900 2545 1715 2420 1900 2678 1642 2480 1900 2793 1587 2540 1900 2814 1626 2600 1900 2835 1665 2660 1900 2856 1704 2720 1900 2877 1743 2780 1900 2898 1782 2840 1900 2918 1822 2900 1900 2939 1861
% 
\special{pn 4}%
\special{pa 2240 1900}%
\special{pa 2278 1864}%
\special{fp}%
\special{pa 2300 1900}%
\special{pa 2412 1790}%
\special{fp}%
\special{pa 2360 1900}%
\special{pa 2546 1716}%
\special{fp}%
\special{pa 2420 1900}%
\special{pa 2678 1642}%
\special{fp}%
\special{pa 2480 1900}%
\special{pa 2794 1588}%
\special{fp}%
\special{pa 2540 1900}%
\special{pa 2814 1626}%
\special{fp}%
\special{pa 2600 1900}%
\special{pa 2836 1666}%
\special{fp}%
\special{pa 2660 1900}%
\special{pa 2856 1704}%
\special{fp}%
\special{pa 2720 1900}%
\special{pa 2878 1744}%
\special{fp}%
\special{pa 2780 1900}%
\special{pa 2898 1782}%
\special{fp}%
\special{pa 2840 1900}%
\special{pa 2918 1822}%
\special{fp}%
\special{pa 2900 1900}%
\special{pa 2940 1862}%
\special{fp}%
\end{picture}%
\hspace{6.5truecm}}

\vspace{0.5truecm}

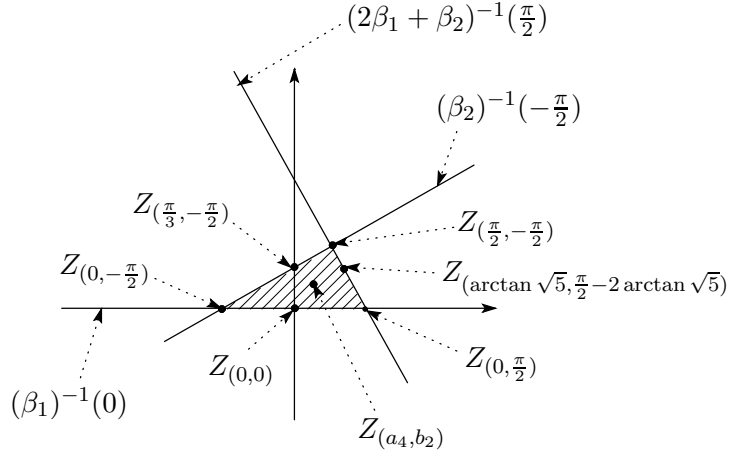
\centerline{\bf Figure 9.}

\vspace{0.5truecm}

$$\begin{tabular}{|c|c|c|c|}
\hline
{\scriptsize$(a,b)$} & {\scriptsize $Z_{(a,b)}$} & 
{\scriptsize$M=(G_2^2)^{\ast}({\rm Exp}\,Z_{(a,b)})$} & 
{\scriptsize ${\rm dim}\,M$}\\
\hline
{\scriptsize $(0,-\frac{\pi}{2})$} & 
{\scriptsize as in Theorem F} & {\scriptsize one-point set} & $0$\\
\hline
{\scriptsize $(0,\frac{\pi}{2})$} & 
{\scriptsize as in Theorem F} & {\scriptsize one-point set} & $0$\\
\hline
{\scriptsize $(\frac{\pi}{2},-\frac{\pi}{2})$} & 
{\scriptsize as in Theorem F} & {\scriptsize totally geodesic} & $8$\\
\hline
{\scriptsize $(\frac{\pi}{3},-\frac{\pi}{2})$} & 
{\scriptsize as in Theorem C} & {\scriptsize not austere} & $6$\\
\hline
{\scriptsize $(\arctan\sqrt 5,\frac{\pi}{2}-2\arctan\sqrt 5)$} & 
{\scriptsize not as in Theorems C$\sim$F} & {\scriptsize not austere} & $10$\\
\hline
{\scriptsize $(a_5,b_3)$} & 
{\scriptsize not as in Theorems C$\sim$F} & {\scriptsize not austere} & $12$\\
\hline
\end{tabular}$$

\vspace{0.2truecm}

\centerline{$(G_2^2)^{\ast}\curvearrowright(G_2\times G_2)/G_2$}

\centerline{$({\rm dim}\,(G_2\times G_2)/G_2=14)$}

\vspace{0.5truecm}

\centerline{{\bf Table 14.}}

\vspace{0.5truecm}

The positions of $Z_0$'s in Table 14 are as in Figure 9.  
The numbers $a_5$ and $b_3$ in Table 14 are real numbers such that 
$a_4,b_2\not\equiv\frac{\pi}{6},\frac{\pi}{3},\frac{\pi}{4},
\frac{3\pi}{4}\,\,\,({\rm mod}\,\pi)$.  

\vspace{1truecm}

\noindent
{\bf Acknowledgement} The author wishes to thank the referee for his valuable comments.  

\vspace{1truecm}

\centerline{{\bf References}}

\vspace{0.5truecm}

{\small 
\noindent
[B] R. Bott, The index theorem for homogeneous differential operators, 
Differential and 

Combinatorial Topology, Princeton University Press, 1965, 167-187.  

\noindent
[BCO] J. Berndt, S. Console and C. Olmos, Submanifolds and Holonomy, 
Research Notes 

in Mathematics 434, CHAPMAN $\&$ HALL/CRC Press, Boca Raton, London, New York 

Washington, 2003.

\noindent
[Co] L. Conlon, Remarks on commuting involutions, Proc. Amer. Math. Soc. 
{\bf 22} (1969) 

255-257.

\noindent
[GT] O. Goertsches and G. Thorbergsson, 
On the Geometry of the orbits of Hermann actions, 

Geom. Dedicata {\bf 129} (2007) 101-118.

\noindent
[He] S. Helgason, 
Differential geometry, Lie groups and Symmetric Spaces, Academic Press, 

New York, 1978.

\noindent
[HPTT] E. Heintze, R.S. Palais, C.L. Terng and G. Thorbergsson, 
Hyperpolar actions on 

symmetric spaces, Geometry, topology and physics for Raoul Bott 
(ed. S. T. Yau), 

Conf. Proc. Lecture Notes Geom. Topology {\bf 4}, Internat. Press, Cambridge, 
MA, 1995 

pp214-245.

\noindent
[HTST] D. Hirohashi, H. Tasaki, H. Song and R. Takagi, Minimal orbits 
of the isotropy 

groups of symmetric spaces of compact type, Differential Goem. Appl. {\bf 13} 
(2000) 

167-177.

\noindent
[HIT] D. Hirohashi, O. Ikawa and H. Tasaki, Orbits of isotropy groups 
of compact symme-

tric spaces, Tokyo J. Math. {\bf 24} (2001) 407-428.

\noindent
[I] O. Ikawa, The geometry of symmetric triad and orbit spaces of Hermann 
actions, 

J. Math. Soc. Japan {\bf 63} (2011) 79-139.  

\noindent
[IST] O. Ikawa, T. Sakai and H. Tasaki, Orbits of Hermann actions, Osaka J. 
Math. {\bf 38} 

(2001) 923-930.  

\noindent
[K1] N. Koike, Actions of Hermann type and proper complex equifocal 
submanifolds, 

Osaka J. Math. {\bf 42} (2005) 599-611.

\noindent
[K2] N. Koike, Collapse of the mean curvature flow for equifocal 
submanifolds, Asian J. 

Math. {\bf 15} (2011) 101-128.
%arXiv:math.DG/0908.3086v5.

\noindent
[Kol] A. Kollross, A Classification of hyperpolar and cohomogeneity one 
actions, Trans. 

Amer. Math. Soc. {\bf 354} (2001) 571-612.

\noindent
[TT] C.L. Terng and G. Thorbergsson, 
Submanifold geometry in symmetric spaces, J. Diff-

erential Geometry {\bf 42} (1995) 665-718.

\vspace{0.5truecm}

{\small 
\rightline{Department of Mathematics, Faculty of Science}
\rightline{Tokyo University of Science, 1-3 Kagurazaka}
\rightline{Shinjuku-ku, Tokyo 162-8601 Japan}
\rightline{(koike@ma.kagu.tus.ac.jp)}
}

\end{document}